\documentclass[reqno]{amsart} 
\usepackage{graphicx}
\usepackage{amsmath}
\usepackage{amssymb}
\usepackage{amsthm}
\usepackage{amstext}
\usepackage{amsbsy}
\usepackage{amsfonts}
\usepackage{enumerate}
\theoremstyle{plain}
\newtheorem{theo}{Theorem}[section]
\theoremstyle{definition}
\theoremstyle{remark}
\numberwithin{equation}{section}

\newcommand{\R}{\mathbb{R}}

\newcommand{\divrg}{\textrm{div}\,}

\newcommand{\onehalf}{\frac{1}{2}}
\newcommand{\body}{\Omega}
\newcommand{\boundary}{\partial \body}

\title[Computing volume bounds of inclusions by EIT measurements] {Computing volume bounds of inclusions \\ by  EIT measurements$^\ast$}

\author[Alessandrini, Bilotta, Morassi, Rosset and Turco] 
{Giovanni~Alessandrini$^\vartriangle$, Antonio~Bilotta$^\circ$,  
Antonino~Morassi$^\triangledown$, Edi~Rosset$^\vartriangle$ and  
Emilio~Turco$^\star$}

\address{$^\vartriangle$ Dipartimento di Matematica e Informatica,  Universit\`a degli Studi di Trieste, Trieste, Italy.}
\address{$^\circ$ Dipartimento di Strutture, Universit\`a della Calabria, Rende (CS), Italy.}
\address{$^\triangledown$ Dipartimento di Georisorse e Territorio, Universit\`a degli Studi di Udine, Udine, Italy.}
\address{$^\star$ Dipartimento di Architettura e Pianificazione, Universit\`a degli Studi di Sassari, Alghero, Italy.}

\subjclass[2000]{35R30, 35R25, 73C02.}

\keywords{Size estimates, electrical impedance tomography}
\date{December 22, 2006}
\thanks{$^\ast$Work supported by MIUR, PRIN no. 2004011204}

\begin{document}

\begin{abstract}
The size estimates approach for Electrical Impedance Tomography (EIT) allows for estimating the size (area or volume) of an unknown inclusion in an electrical conductor by means of one pair of boundary measurements of voltage and current. In this paper we show by numerical simulations how to obtain such bounds for practical application of the method. The computations are carried out both in a 2--D and a 3--D setting.
\end{abstract}
\maketitle
\section{Introduction} \label{sec:introduction}

EIT is aimed at imaging the internal conductivity of a body from
current and voltage measurements taken at the boundary. It is well
known, \cite{l:a88}, \cite{l:m01}, that, even in the ideal
situation in which all possible boundary measurements are
available, the correspondence \emph{boundary data} $\rightarrow$
\emph{conductivity} is highly (exponentially) unstable. As a
consequence it is evident that, in practice, it is impossible to
distinguish high resolution features of the interior from limited
and noisy boundary data, \cite{l:av}.

Motivated by applications, a line of investigation pursued by many
authors, \cite{l:fr}, \cite{l:frg}, \cite{l:fri}, \cite{l:ai},
\cite{l:fks}, \cite{l:aip}, \cite{l:isak}, \cite{l:isak-libro}, has been the one of limiting the analysis to cases
in which one seeks an unknown interior inclusion embedded in an
otherwise known (may be even homogeneous) conductor, and whose
conductivity is assumed to differ from the background.

Even in this restricted case, and even when full boundary data are
available, the instability remains of exponential type
\cite{l:dcr}.

It is therefore reasonable to further restrict the goal and attempt to evaluate some parameters expressing the size (area,
volume) of the inclusion, disregarding its precise location and shape, having at our disposal one pair of boundary measurements of voltage and current. This approach, which can be traced back to \cite{l:fr}, has been well developed theoretically, \cite{l:ar98},
\cite{l:kss}, \cite{l:ars}, \cite{l:amr03}, see also \cite{l:ike} and \cite{l:amr04} for the analogous treatment in the linear elasticity framework. In order to describe such type of results we need first to introduce some notation.

We denote by $\Omega$ a bounded domain in $\R^n$, $n=2,3$,
representing an electrical conductor. The boundary $\partial
\Omega$ of $\Omega$ is assumed of Lipschitz class, with constants
$r_0$, $M_0$, that is the boundary can be locally represented as a
graph of a Lipschitz continuous function with Lipschitz constant
$M_0$ in some ball of radius $r_0$. When no inclusion is present
in the conductor we assume that it is homogeneous and we pose its
conductivity $\sigma(x)\equiv 1$. When the conductor contains an
unknown inclusion $D$ of different conductivity, say $k>0$, $k
\neq 1$ the overall conductivity in the conductor will be given by
$\sigma(x)=1+(k-1)\chi_D(x)$. Here and in what follows it is
assumed that $D$ is strictly contained in $\Omega$. More
precisely, for a given $d_0> 0$,
\begin{equation}
   \label{eq:2.condition_d0}
   \textrm{dist}(D, \partial \Omega) \geq d_0.
\end{equation}
Let $\varphi \in H^{- \frac{1}{2}}(\partial \Omega)$,
$\int_{\partial \Omega} \varphi =0$, be an applied current density
on $\partial \Omega$. The induced electrostatic potential $u \in
H^1(\Omega)$ is the solution of the Neumann problem
\begin{equation}
   \label{eq:2.Neumann_pbm_with_incl}
   \left\{ \begin{array}{ll}
   \divrg ((1+(k-1) \chi_D) \nabla u)=0, &
   \mathrm{in}\ \Omega ,\\
    &  \\
   \nabla u \cdot \nu= \varphi, &
   \mathrm{on}\ \partial \Omega,
   \end{array}\right.
\end{equation}
where $\nu$ denotes the outer unit normal to $\partial \Omega$.

When $D$ is the empty set, that is when the inclusion is absent,
the reference electrostatic potential $u_0 \in H^1(\Omega)$
satisfies the Neumann problem
\begin{equation}
   \label{eq:2.Neumann_pbm_without_incl}
   \left\{ \begin{array}{ll}
   \Delta u_0=0, &
   \mathrm{in}\ \Omega ,\\
    &  \\
   \nabla u_0 \cdot \nu= \varphi, &
   \mathrm{on}\ \partial \Omega.
   \end{array}\right.
\end{equation}

In both cases \eqref{eq:2.Neumann_pbm_with_incl} and
\eqref{eq:2.Neumann_pbm_without_incl}, the solutions $u$ and $u_0$
are determined up to an additive constant.

Let us denote by $W$, $W_0$ the powers required to maintain the
current density $\varphi$ on $\partial \Omega$ when $D$ is present
or it is absent, respectively. Namely
\begin{equation}
   \label{eq:2.def_W}
   W=\int_{\partial \Omega} u \varphi = \int_{\Omega}(1+(k-1)\chi_D)| 
\nabla u|^2,
\end{equation}
\begin{equation}
   \label{eq:2.def_W0}
   W_0=\int_{\partial \Omega} u_0 \varphi = \int_{\Omega}|\nabla u_0|^2.
\end{equation}
The size estimate approach developed in \cite{l:ar98},
\cite{l:kss}, \cite{l:ars}, \cite{l:amr03}, tells us that the
measure $|D|$ of $D$ can be bounded from above and below in
terms of the quantity $\left|\frac{W_0-W}{W_0}\right|$ which we
call the normalized power gap. More precisely the following bounds
hold, see \cite[Theorem 2.3]{l:amr03}.

\begin{theo}
   \label{theo:size-estim-EIT-general}
   Let $D$ be any measurable subset of $\Omega$ satisfying \eqref{eq:2.condition_d0}. Under the above assumptions, if $k > 1$ we have
\begin{equation}
   \label{eq:2.size-estim-EIT-more-conduct}
   \frac {1} {k-1} C^{+}_{1}
   \frac{W_0-W}{W_0}
   \leq
   |D|
   \leq
   \left (
   \frac{k}{k-1}
   \right )^{ \frac{1}{p} }
   C^{+}_{2}
   \left (
   \frac{W_0-W}{W_0}
   \right )^{ \frac{1}{p} }.
   \end{equation}
If, conversely, $k < 1$, then we have
\begin{equation}
   \label{eq:2.size-estim-EIT-less-conduct}
   \frac {k} {1-k} C^{-}_{1}
   \frac{W-W_0}{W_0}
   \leq
   |D|
   \leq
   \left (
   \frac{1}{1-k}
   \right )^{ \frac{1}{p} }
   C^{-}_{2}
   \left (
   \frac{W-W_0}{W_0}
   \right )^{ \frac{1}{p} },
   \end{equation}
where $C^{+}_{1}$, $C^{-}_{1}$ only depend on $d_0$, $|\Omega|$,
$r_0$, $M_0$, whereas $p>1$, $C^{+}_{2}$, $C^{-}_{2}$ only depend
on the same quantities and, in addition, on the \textit{frequency
of $\varphi$}
\begin{equation}
   \label{eq:2.frequency}
   F[\varphi] = \frac{\|\varphi \|_{H^{ -\frac{1}{2} }(\partial  
\Omega)}}{\|\varphi \|_{H^{-1}(\partial \Omega)}}.
\end{equation}
\end{theo}

When it is a priori known that the inclusion $D$ is not too small
(if it is at all present), a situation which often occurs in
practical applications, stronger bounds apply.

\begin{theo}
   \label{theo:size-estim-EIT-fat-incl}
   Under the above hypotheses, let us assume, in addition, that
\begin{equation}
   \label{eq:2.fat-inclusion}
   |D| \geq m_0,
\end{equation}
for a given positive constant $m_0$. If $k > 1$ we have
\begin{equation}
   \label{eq:2.size-estim-EIT-more-conduct-fat-incl}
   \frac {1} {k-1} C^{+}_{1}
   \frac{W_0-W}{W_0}
   \leq
   |D|
   \leq
   \frac{k}{k-1}
   C^{+}_{2}
   \frac{W_0-W}{W_0}
   .
   \end{equation}
If, conversely, $k < 1$, then we have
\begin{equation}
   \label{eq:2.size-estim-EIT-less-conduct-fat-incl}
   \frac {k} {1-k} C^{-}_{1}
   \frac{W-W_0}{W_0}
   \leq
   |D|
   \leq
   \frac{1}{1-k}
   C^{-}_{2}
   \frac{W-W_0}{W_0},
   \end{equation}
where $C^{+}_{1}$, $C^{-}_{1}$ only depend on $d_0$, $|\Omega|$,
$r_0$, $M_0$, whereas $C^{+}_{2}$, $C^{-}_{2}$ only depend on the
same quantities and, in addition, on $m_0$ and $F[\varphi]$.
\end{theo}

Theorem \ref{theo:size-estim-EIT-fat-incl} can be easily deduced
from Theorem \ref{theo:size-estim-EIT-general} by the arguments
sketched in \cite[Appendix]{l:abfmrt04}.

One of the goals of the present paper is to test the applicability
of such bounds by numerical simulations with the following
purposes:

\emph{i) provide practical evaluations of the constants
$C_1^{\pm}$, $C_2^{\pm}$ appearing in the above inequalities
\eqref{eq:2.size-estim-EIT-more-conduct},
\eqref{eq:2.size-estim-EIT-less-conduct},
\eqref{eq:2.size-estim-EIT-more-conduct-fat-incl},
\eqref{eq:2.size-estim-EIT-less-conduct-fat-incl}; }

\emph{ii) when, due to special geometric configurations, it is
possible to compute theoretically such constants, compare such
theoretical values with those obtained by simulations;}

\emph{iii) show that such upper and lower bounds deteriorate as
the frequency $F[\varphi]$ increases.}

The other goal of this paper is to perform similar kinds of
numerical simulations when the so-called \emph{complete model} of
EIT is adopted. We recall that this model is aimed at an accurate
description of the boundary measurements suitable for medical
applications, and was introduced in \cite{l:cing} and subsequently
developed in \cite{l:pbp} and \cite{l:sci}. In this model, the
metal electrodes behave as perfect conductors and provide a
low-resistance path for current. An electrochemical effect at the
contact between the electrodes and the body results in a thin,
highly resistive, layer. The impedance of this layer is
characterized by a positive quantity $z_l$ on each electrode
$e_l$, $l=1,...,L$, which is called \textit{surface impedance}.

Denoting by $I_l$ the current applied to each $e_l$, the resulting
boundary condition on each electrode $e_l$ becomes
\begin{equation}
   \label{eq:2.bound-cond}
   u+z_l \nabla u \cdot \nu = U^l, \quad \quad \hbox{on } e_l,
\end{equation}
where the unknown constant $U^l$ is the voltage which can be
measured at the electrode $e_l$.

We assume, as before, that the reference conductor has
conductivity $\sigma \equiv 1$ and that an unknown inclusion $D$
of conductivity $\sigma \equiv k$, with $k>0$ and $k \neq 1$, is
strictly contained in $\Omega$. Therefore, the electrostatic
potential $u$ inside the conductor is determined, up to an
additive constant, as the solution to the following problem
\begin{equation}
   \label{eq:2.Phys-Neumann_pbm_with_incl}
   \left\{ \begin{array}{ll}
   \divrg ((1+(k-1) \chi_D) \nabla u)=0, &
   \mathrm{in}\ \Omega ,\\
    u+z_l \nabla u \cdot \nu = U^l, &
   \mathrm{on}\ e_l, \ 1 \leq l \leq L, \\
   \nabla u \cdot \nu= 0, &
   \mathrm{on}\ \partial \Omega \setminus \cup_{l=1}^L e_l, \\
   \int_{e_l} \nabla u \cdot \nu = I_l, &
   \ 1 \leq l \leq L,
   \end{array}\right.
\end{equation}
where the so-called current pattern $I=(I_1, ..., I_L)$ is subject
to the conservation of charge condition $\sum_{l=1}^L I_l=0$, and
the unknown constants $U^l$ are the components of the so-called
voltage pattern $U=(U^1, ..., U^L)$.

When the inclusion is absent, the electrostatic potential $u_0$
induced by the same current pattern $I$ is determined, up to an
additive constant, as the solution of the following problem
\begin{equation}
   \label{eq:2.Phys-Neumann_pbm_without_incl}
   \left\{ \begin{array}{ll}
   \Delta u_0=0, &
   \mathrm{in}\ \Omega ,\\
    u_0+z_l \nabla u_0 \cdot \nu = U_0^l, &
   \mathrm{on}\ e_l, \ 1 \leq l \leq L, \\
   \nabla u_0 \cdot \nu= 0, &
   \mathrm{on}\ \partial \Omega \setminus \cup_{l=1}^L e_l, \\
   \int_{e_l} \nabla u_0 \cdot \nu = I_l, &
   \ 1 \leq l \leq L,
   \end{array}\right.
\end{equation}
where, as before, the $U_0^l$ are unknown constants in the direct
problem \eqref{eq:2.Phys-Neumann_pbm_without_incl}.

We shall assume that the sets $e_1,...,e_L$, representing the
electrodes, are open, pairwise disjoint, connected subsets of
$\partial \Omega$ and, in addition,
\begin{equation}
   \label{eq:2.cond-electrodes}
   \textrm{dist}(e_l,e_k) \geq \delta_1 > 0 \quad \hbox{for every }  
l,k, \ l \neq k.
\end{equation}

The surface impedance $z_l$ on $e_l$, $l=1,...,L$, is assumed to
be real valued and to satisfy the following bounds
\begin{equation}
   \label{eq:2.cond-impedance}
   0<m\leq z_l\leq M,  \quad \hbox{for every } l=1,...,L.
\end{equation}

In this formulation, the powers $W$ and $W_0$ become
\begin{equation}
   \label{eq:2.def_W-phys-model}
   W=\sum_{i=1}^L I_i U^i,
\end{equation}
\begin{equation}
   \label{eq:2.def_W0-phys-model}
   W_0=\sum_{i=1}^L I_i U_0^i.
\end{equation}
Size estimates like those of Theorems
\ref{theo:size-estim-EIT-general},
\ref{theo:size-estim-EIT-fat-incl} were obtained for the complete
model in \cite{l:ar04}. In particular we have

\begin{theo}
   \label{theo:size-estim-Phys-EIT-general}
   Let $D$ be any measurable subset of $\Omega$ satisfying \eqref{eq:2.condition_d0} and let $W$, $W_0$ be given by \eqref{eq:2.def_W-phys-model},  \eqref{eq:2.def_W0-phys-model}. Then, inequalities \eqref{eq:2.size-estim-EIT-more-conduct}, \eqref{eq:2.size-estim-EIT-less-conduct} hold for  $k > 1$ and $k < 1$, respectively, where the constants $C^{+}_{1} $, $C^{-}_{1}$ only depend on $d_0$, $|\Omega|$, $r_0$,  $M_0$, and $C^{+}_{2}$, $C^{-}_{2}$ and $p>1$ only depend on the same quantities and, in addition,  on $\delta_1$, $M$ and $m$.
\end{theo}

Also in this case, the size estimates of $|D|$ can be improved to the form \eqref{eq:2.size-estim-EIT-more-conduct-fat-incl},
\eqref{eq:2.size-estim-EIT-less-conduct-fat-incl} when condition
\eqref{eq:2.fat-inclusion} is satisfied.

\medskip

In Section \ref{sec:num-EIT} we consider the standard EIT setting.
We start by describing the finite element setup used in our
numerical simulations in Section \ref{subsec:nummodel}. Next (as a
warmup) we consider a two-dimensional model in Section
\ref{subsec:2D}.

In Section \ref{subsec:3D} we consider the three-dimensional case
and we discuss all items i), ii), iii) introduced above. In
particular we observe that, comparing the results as the frequency
$F[\varphi]$ increases, we have quite rapidly a serious
deterioration of the bounds. This poses a severe warning on the
limitations that have to be taken into account in the choice of
the boundary current profile $\varphi$.

Section \ref{sec: num-phys-EIT} is devoted to simulations with the
complete EIT model. In this case it is reasonable analyze the case
when only two electrodes, one positive and one negative, are
attached to the surface of the conductor. In this case, the
frequency function is not available {}from the data since we are
not prescribing the boundary current $\nabla
u\cdot\nu_{|\partial\Omega}$ but only the current pattern, which
is a 2-electrode configuration, is just the pair $(1,-1)$. In
place of the frequency function, the parameters that may influence
the constants in the volume bounds are: the width of the
electrodes and the distance between them. We perform various
experiments to test such variability.

\section{Numerical simulations for the EIT model} \label{sec:num-EIT}

\subsection{Numerical model}  \label{subsec:nummodel}

The numerical model is based on the discretization of the energy
functional $J: H^1(\Omega, \R^n) \rightarrow \R$
\begin{equation}
     \label{eq:energy_cont}
     J(u) = \frac{1}{2} \int_\Omega (1 + (k-1) \chi_D ) \nabla u  
\cdot \nabla u - \int_{\partial \Omega}  \varphi u ,
\end{equation}
associated to the variational formulation of problem
(\ref{eq:2.Neumann_pbm_with_incl}). The energy functional
(\ref{eq:energy_cont}) has been discretized by using the High
Continuity (HC) technique already presented in \cite{l:ari85} and
\cite{l:bft04} in the context of linear elasticity. Accordingly,
for 2--D problems the electric potential on the $e$--th
finite element can be represented as
\begin{equation}
     \label{eq:HCdispl}
     u_e = \sum_{i,j=1}^{3} \phi_i(\xi_1) \phi_j(\xi_2) u_{ij},
\end{equation}
whereas for the 3--D case it assumes the form
\begin{equation}
     \label{eq:HCdispl3D}
     u_e = \sum_{i,j,l=1}^{3} \phi_i(\xi_1) \phi_j(\xi_2) \phi_l 
(\xi_3) u_{ijl},
\end{equation}
where the coordinates $\xi_r$, $r=1,...,n$, span the unitary
element domain $[-\frac{1}{2},\frac{1}{2}]^n$, $n=2,3$, and
$u_{ij}$, $u_{ijl}$ are the HC parameters involved in the field
interpolation on the generic element. The shape functions
$\phi_i(\xi_r)$ are defined as
\begin{equation} \label{eq:HCfun}
     \left\{
     \begin{split}
     \phi_1(\xi_r)=& \frac{1}{8} - \onehalf \xi_r + \onehalf \xi_r^2, \\
     \phi_2(\xi_r)=& \frac{3}{4} - \xi_r^2, \\
     \phi_3(\xi_r)=& \frac{1}{8} + \onehalf \xi_r + \onehalf \xi_r^2.
     \end{split}
     \right.
\end{equation}
\begin{figure}[ht]
     \centering
     \includegraphics[width=10cm]{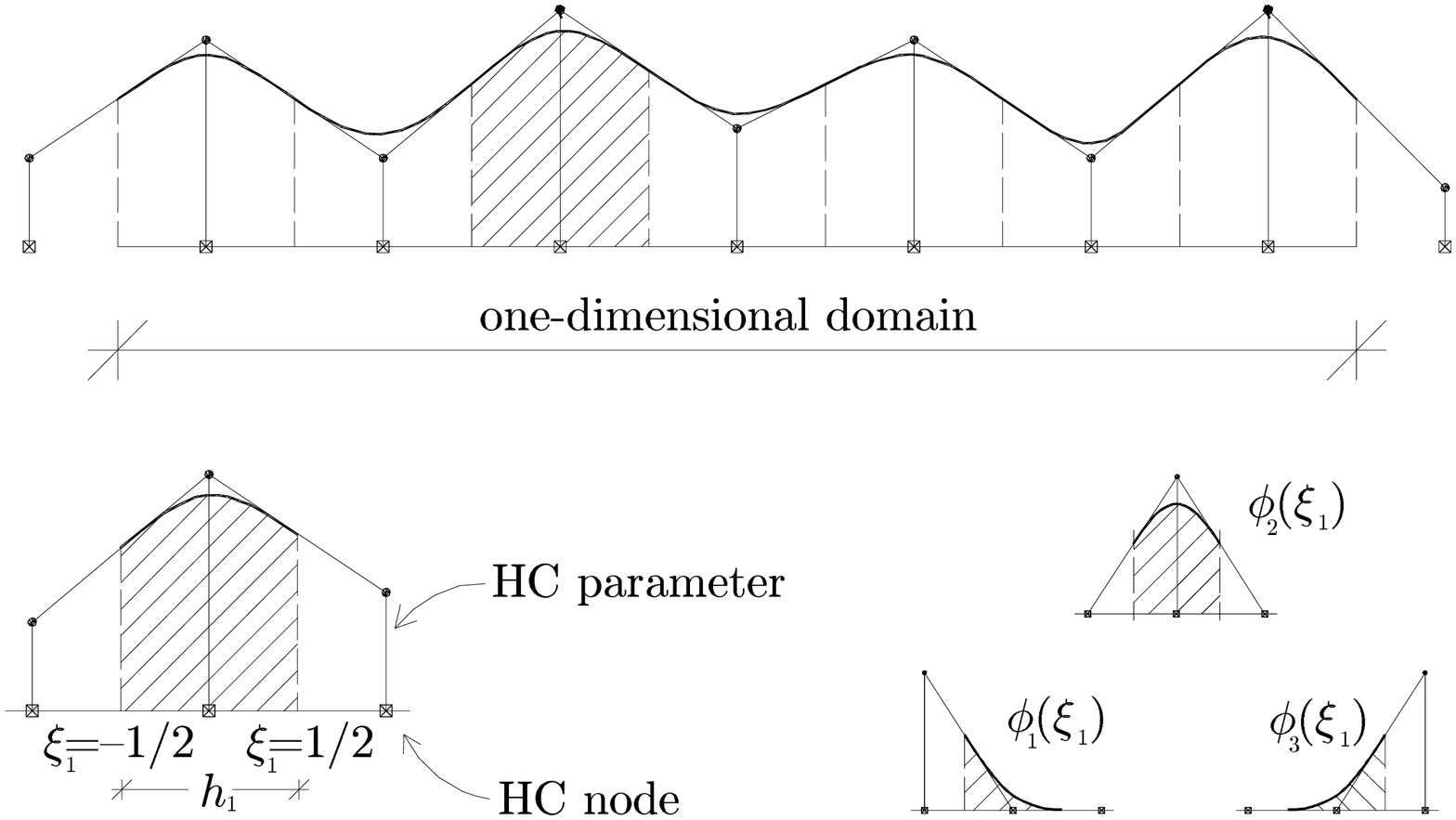} 
     \caption{HC interpolation in the 1--D case: nodes, parameters and shape functions.}
     \label{fig:HCmono}
\end{figure}
\begin{figure}[ht]
     \centering
     \includegraphics[width=12cm]{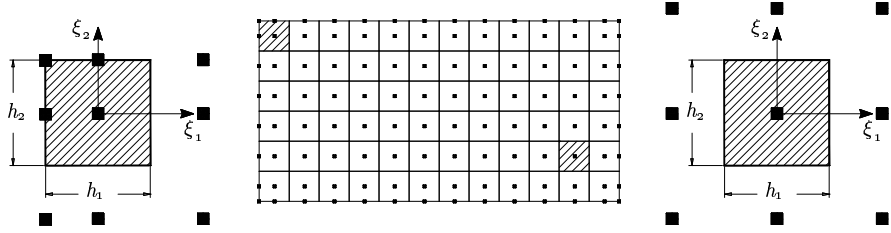}
     \caption{HC mesh in the 2--D case: nodes for boundary and inner
elements.}
     \label{fig:HCbi}
\end{figure}
%
%

The 1--D case illustrated in Figure~\ref{fig:HCmono} shows the
meaning of the HC parameters. They allow to define the slopes of
the interpolated function at the end points of the element. On the
same figure one can see also the positions of the HC nodes and the
shape functions (\ref{eq:HCfun}).

Figure~\ref{fig:HCbi} shows a typical structured mesh on a
rectangular domain and the nodes used for the approximation of the
potential field in the 2--D case. For elements with a side lying
on the boundary, in order to easily impose the Neumann boundary
conditions, special shape functions are used. In practice, the
external HC nodes are translated onto the boundary $\partial
\Omega$ and the related HC parameters have the meaning of function
values (see again Figure~\ref{fig:HCbi}). In this case the shape
functions relative to a {\em left} boundary ($\xi_r=-\frac{1}{2}$)
and a {\em right} boundary ($\xi_r=\frac{1}{2}$) of the finite
element are
\begin{equation}
     \label{eqn:HCfunLR}
     {\rm left:} \left\{
     \begin{split}
     \phi_1(\xi_r)\;\; =& \frac{1}{4} - \xi_r + \xi_r^2, \\
     \phi_2(\xi_r)\;\; =& \frac{5}{8} + \onehalf \xi_r - \frac{3}{2}  
\xi_r^2, \\
     \phi_3(\xi_r)\;\; =& \frac{1}{8} + \onehalf \xi_r + \onehalf
     \xi_r^2;
     \end{split}
     \right. \qquad\; {\rm right:} \left\{
     \begin{split}
     \phi_1(\xi_r)\;\; =& \frac{1}{8} - \onehalf \xi_r + \onehalf  
\xi_r^2, \\
     \phi_2(\xi_r)\;\; =& \frac{5}{8} - \onehalf \xi_r - \frac{3}{2}  
\xi_r^2, \\
     \phi_3(\xi_r)\;\; =& \frac{1}{4} + \xi_r + \xi_r^2.
     \end{split}
     \right.
\end{equation}

Further details about the HC interpolation can be found in
\cite{l:ari85} and \cite{l:bft04}. This interpolation technique,
which can be considered as a particular case of the B\'{e}zier
interpolation, has the main advantage of reproducing potential
fields of $C^1$ smoothness with a computational cost equivalent to
a $C^0$ interpolation.

By \eqref{eq:HCdispl} or \eqref{eq:HCdispl3D}, the potential field
$u$ on each element $e$ takes the compact form
\begin{equation}\label{eq:HCinterp_mat}
     u_e = \mathbf{N}_e \mathbf{w}_e.
\end{equation}
The one--row matrix $\mathbf{N}_e$ collects the shape functions of
the HC interpolation, whereas the components of the vector
$\mathbf{w}_e$ are the nodal parameters of the underlying element.
With this notation, the gradient of the potential field is given
by
\begin{equation}\label{eq:HCinterp_grad}
     \nabla u_e = \nabla \mathbf{N}_e \mathbf{w}_e.
\end{equation}

We remark that the dimensions of the matrices $\mathbf{N}_e$,
$\nabla \mathbf{N}_e$ and vector $\mathbf{w}_e$ are $1 \times 9$,
$2 \times 9$ and $9 \times 1$ for the 2--D case and $1 \times 27$,
$3\times 27$ and $27\times1$ for the 3--D case.

By \eqref{eq:HCinterp_mat} and \eqref{eq:HCinterp_grad}, the
discrete form of \eqref{eq:energy_cont} becomes
\begin{equation}\label{eq:EIT_var_form_discr}
     J( \mathbf{w}_e ) = \sum_e \left( \frac{1}{2}  \int_{\body_e} (1  
+ (k-1) \chi_D) (\nabla \mathbf{N}_e \mathbf{w}_e)
     \cdot (\nabla \mathbf{N}_e \mathbf{w}_e) - \int_{\boundary_e}  
\varphi \mathbf{N}_e \mathbf{w}_e \right) ,
\end{equation}
or, in a compact form,
\begin{equation}\label{eq:EIT_min_w}
     J(\mathbf{w}_e) = \sum_e \left(  \mathbf{w}_e^T \mathbf{K}_e  
\mathbf{w}_e  - \mathbf{w}_e^T \mathbf{p}_e \right).
\end{equation}
The latter equation provides the definition of the matrix and vector associated to $e$--th element
\begin{equation}\label{eq:EIT_mat}
     \left\{
     \begin{split}
       \mathbf{K}_e = & \int_{\body_e} (1 + (k-1) \chi_D) (\nabla  
\mathbf{N}_e)^T \nabla \mathbf{N}_e, \\
       \mathbf{p}_e = & \int_{\boundary_e} \varphi \mathbf{N}_e,
     \end{split}
     \right.
\end{equation}
which can be used to assemble, by using standard techniques, the
system of equations to solve.

\subsection{Two--dimensional case} \label{subsec:2D}

\begin{figure}[ht]
     \begin{minipage}{.48\textwidth}
     \centering
     \includegraphics[width=5.5cm]{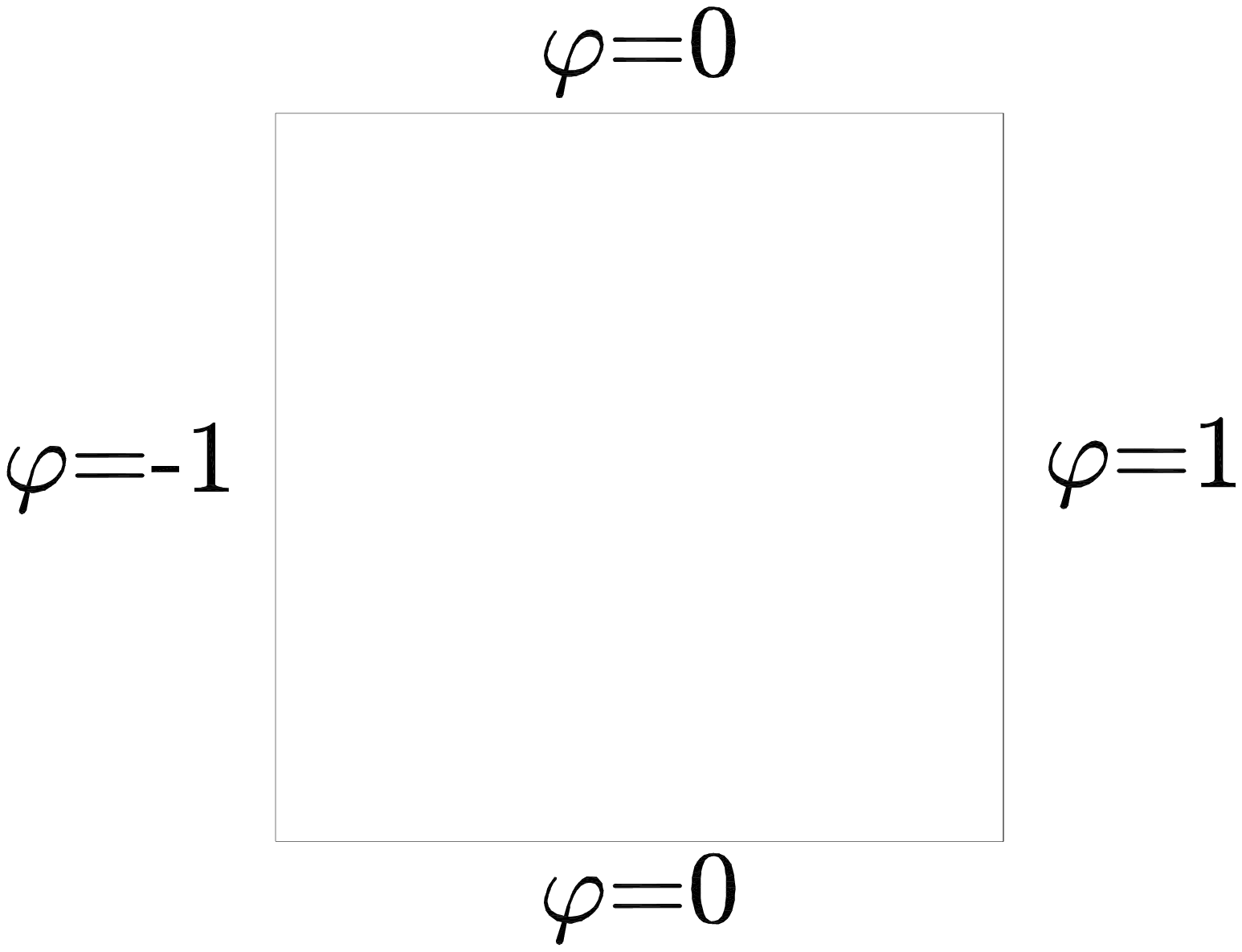}\\
     \centering{(a)} 
   \end{minipage}
   \begin{minipage}{.48\textwidth}
     \centering
     \includegraphics[width=5.5cm]{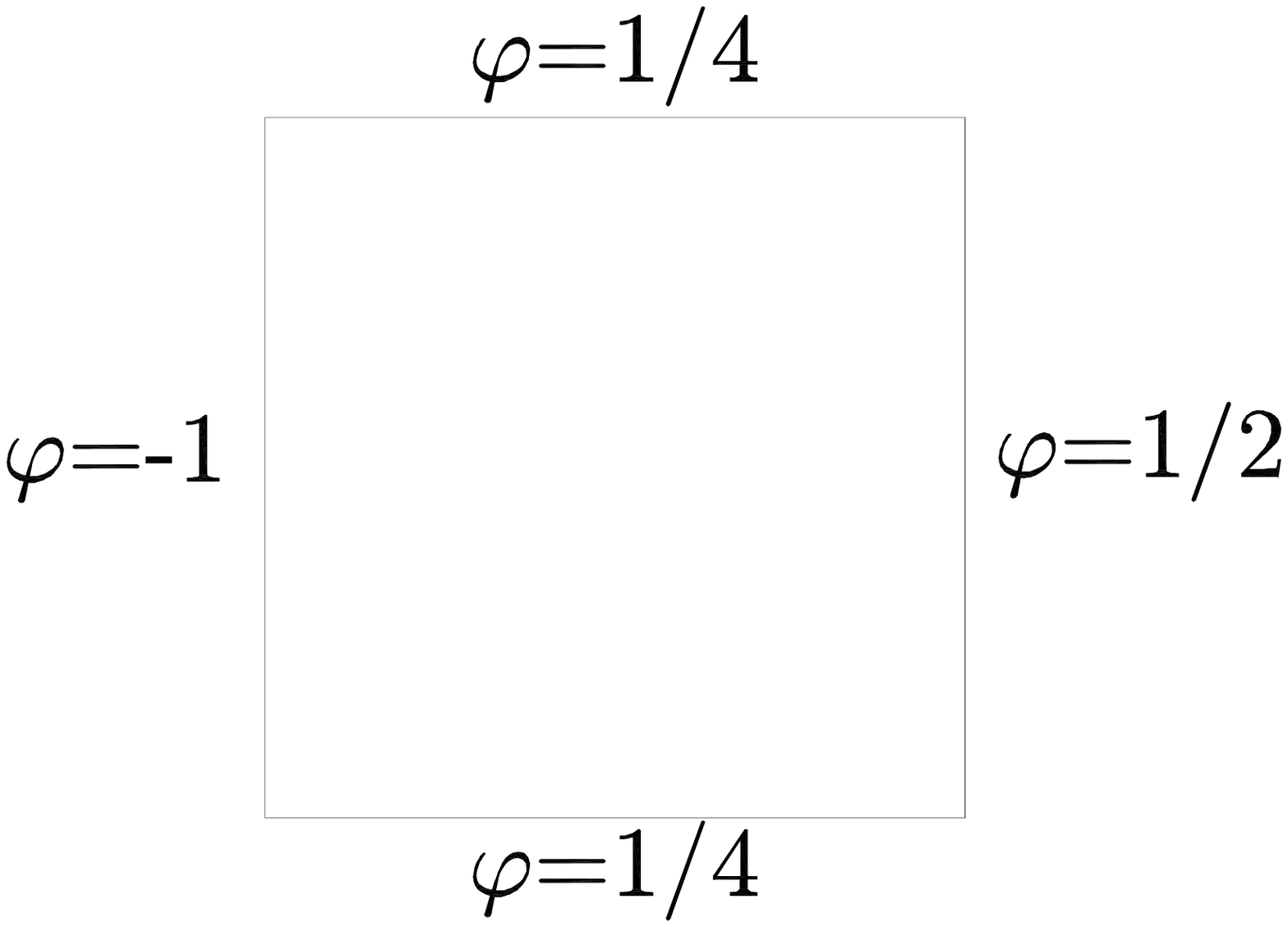}\\
     \centering{(b)} 
   \end{minipage}
   \caption{Square conductor considered in 2--D numerical simulations  
for the EIT model and applied current density fields: Test $T_1$ (a),  
Test $T_2$ (b).}
   \label{fig:Test numerici_2D}
\end{figure}

Numerical analysis has been performed on a square conductor
$\Omega$ of side $l$ under the two current density fields
$\varphi$ illustrated in Figure \ref{fig:Test numerici_2D}. The
domain $\Omega$ has been discretized with a mesh of $21 \times 21$
HC finite elements and for both Test $T_1$ and Test $T_2$ of
Figure \ref{fig:Test numerici_2D} we have considered an inclusion
$D$ with conductivity $k=0.1$ or $k=10$.

\begin{figure}[ht]
   \begin{minipage}{.48\textwidth}
     \centering
     \includegraphics[width=6cm]{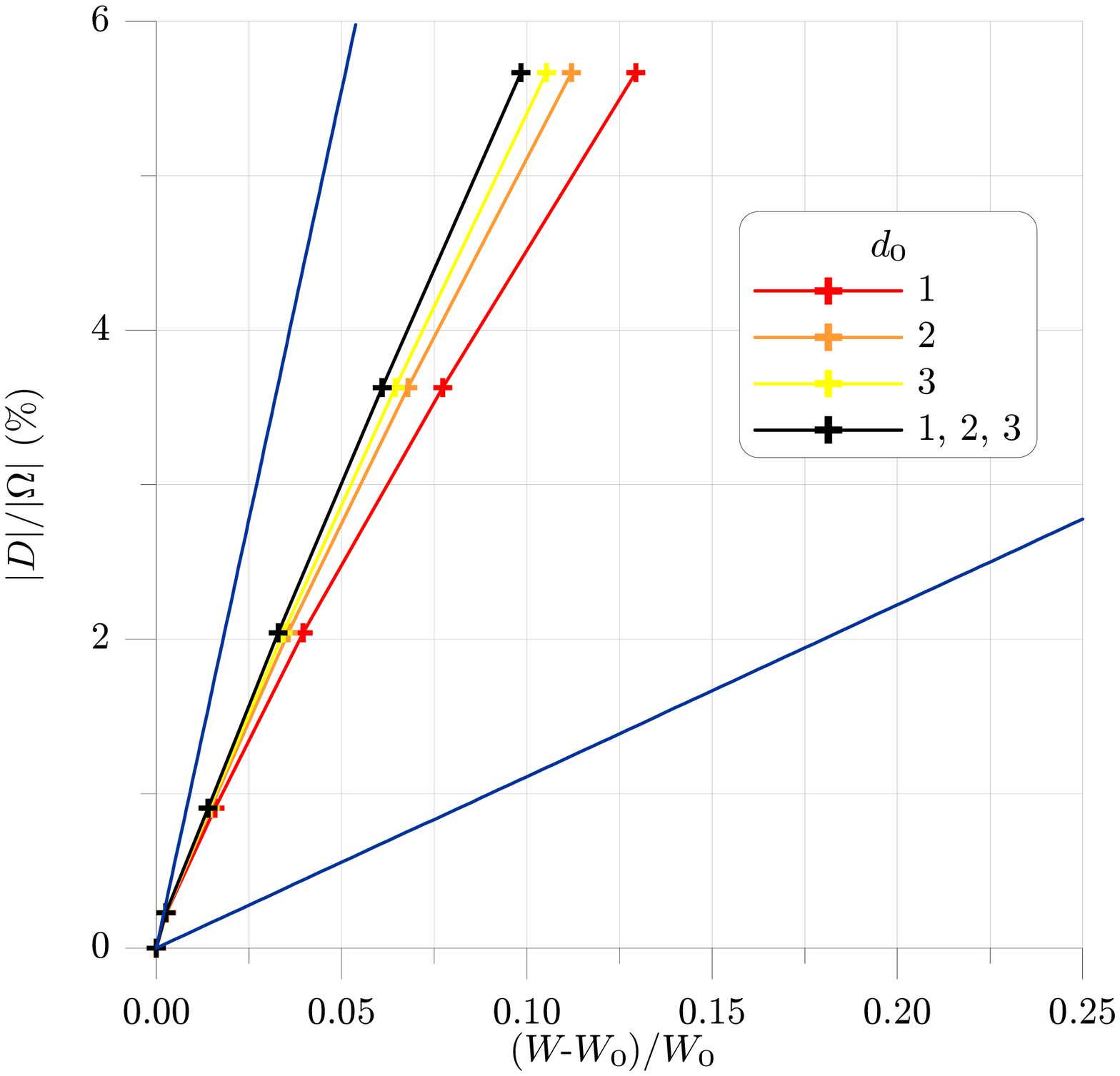}\\
     \centering{(a)} 
   \end{minipage}
   \begin{minipage}{.48\textwidth}
     \centering
     \includegraphics[width=6cm]{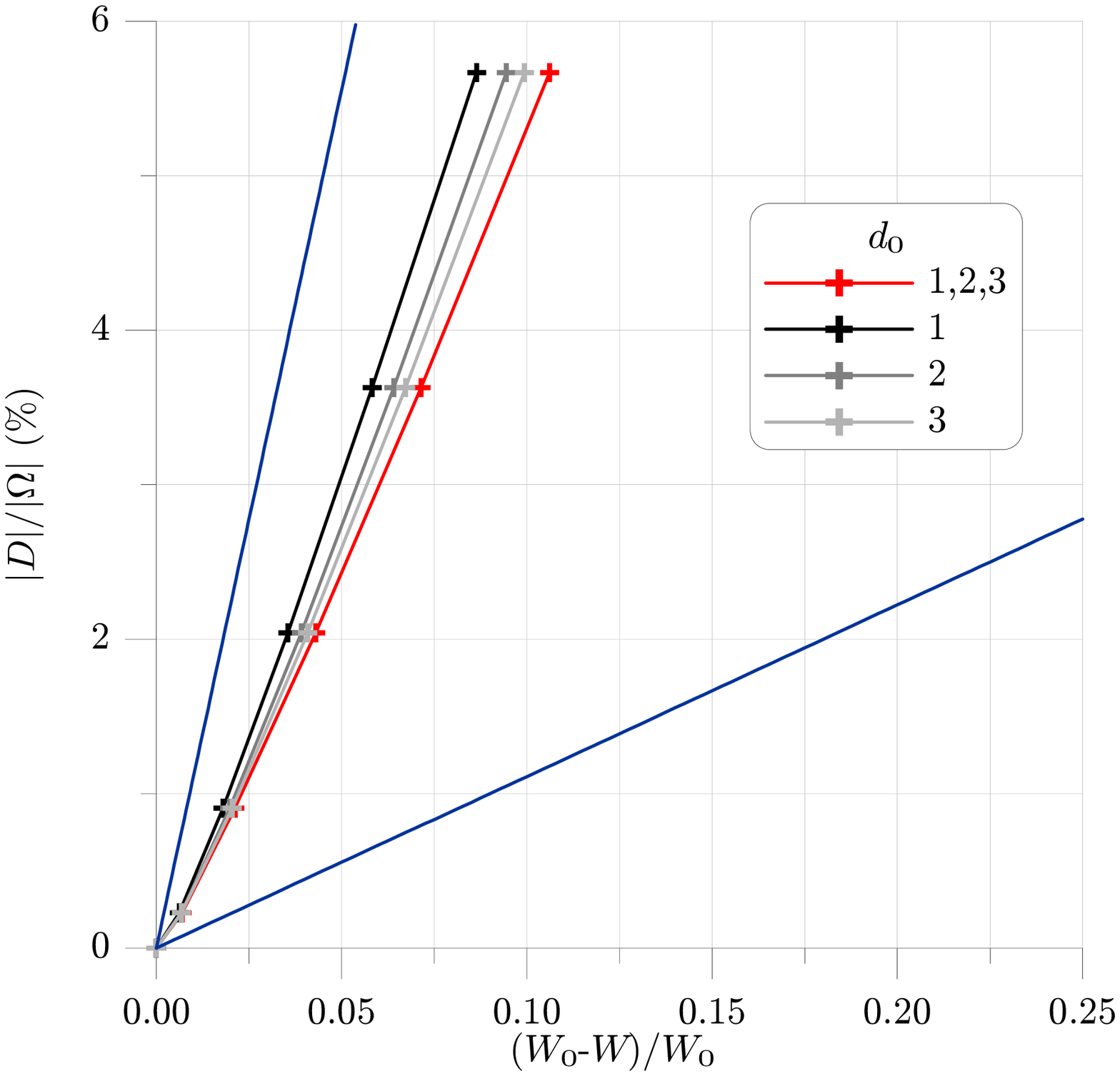}\\
     \centering{(b)} 
   \end{minipage}
   \caption{Influence of $d_0$ for square inclusions in Test $T_1$ of  
Figure \ref{fig:Test numerici_2D}(a) ($21 \times 21$ FE mesh): $k=0.1 
$ (a), $k=10$ (b).}
   \label{fig:T1_2D_pos}
\end{figure}

\begin{figure}[ht]
   \begin{minipage}{.48\textwidth}
     \centering
     \includegraphics[width=6cm]{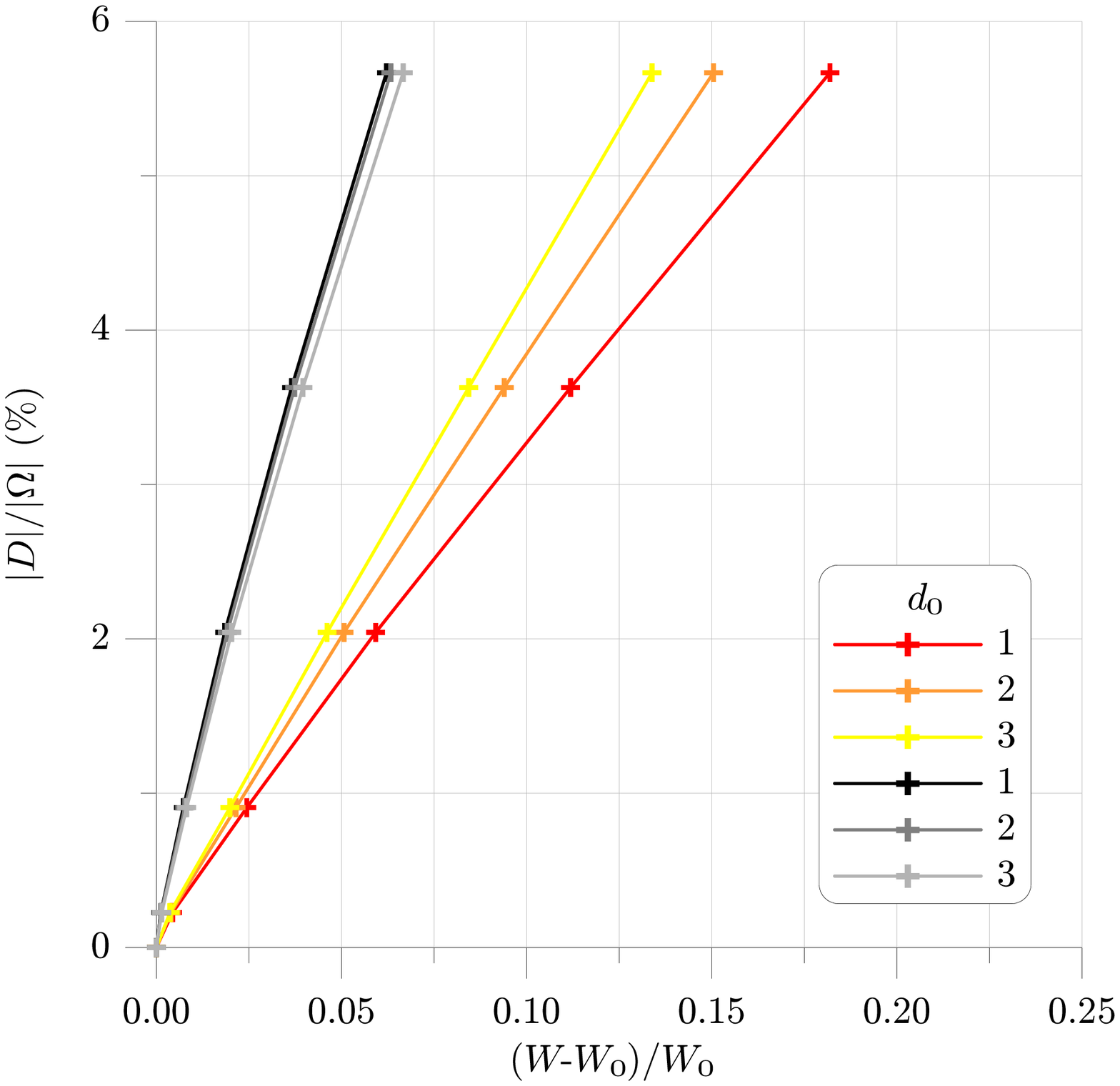}\\
     \centering{(a)}
   \end{minipage}
   \begin{minipage}{.48\textwidth}
     \centering
     \includegraphics[width=6cm]{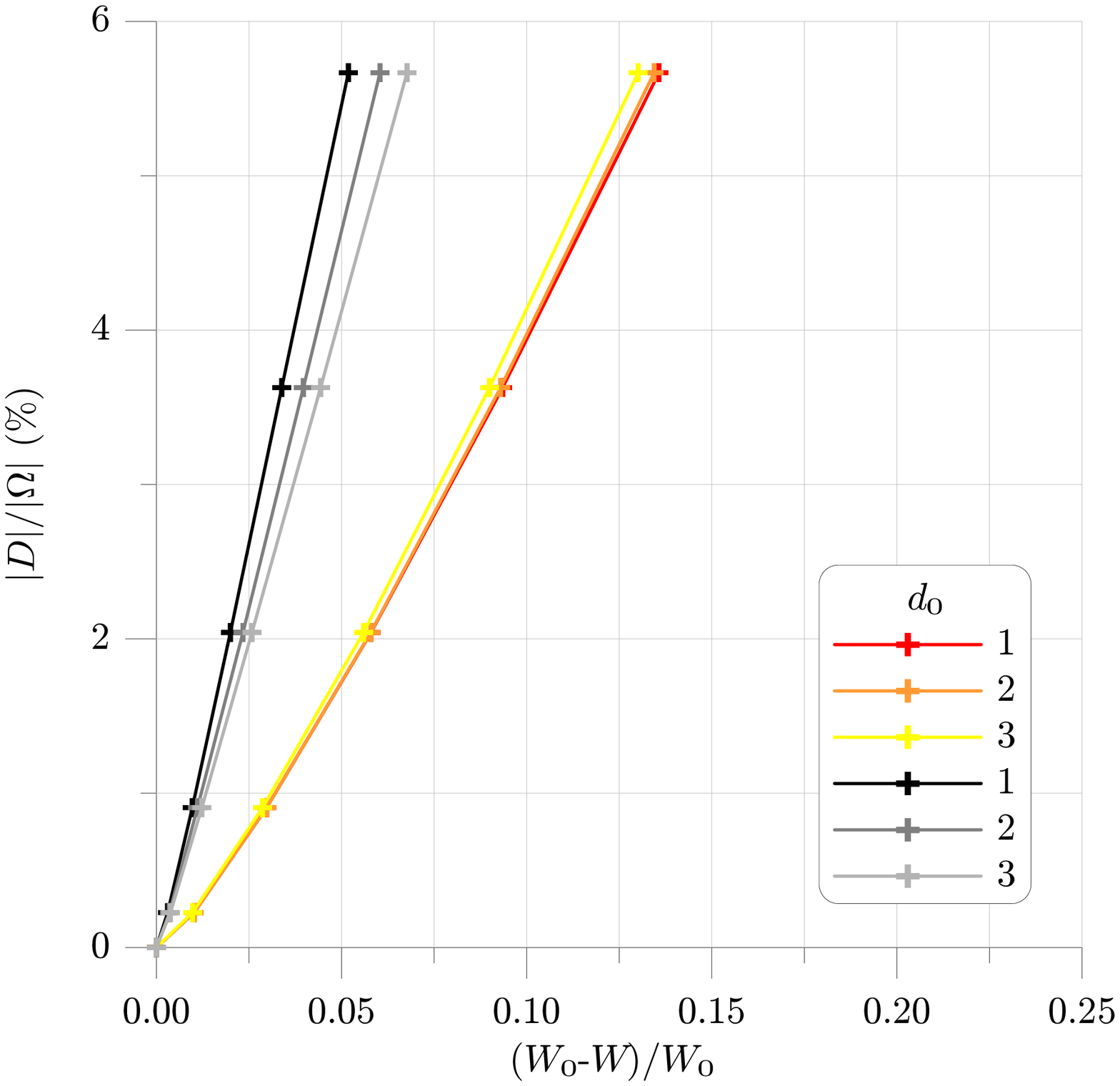}\\
     \centering{(b)}
   \end{minipage}
   \caption{Influence of $d_0$ for square inclusions in Test $T_2$ of  
Figure \ref{fig:Test numerici_2D}(b) ($21 \times 21$ FE mesh): $k=0.1 
$ (a), $k=10$ (b).}
   \label{fig:T2_2D_pos}
\end{figure}

A first series of experiments has been carried out by considering
all possible square inclusions with side ranging from $1$ to $5$
elements, that is the size of inclusion has been kept lower than
$6\%$ of the total size of the conductor. The results are
collected in Figures \ref{fig:T1_2D_pos} and \ref{fig:T2_2D_pos}
for different values of the minimum distance $d_0$ between the
inclusion $D$ and the boundary of $\Omega$.


From Figures \ref{fig:T1_2D_pos}(a) and \ref{fig:T2_2D_pos}(a),
which refer to the case $k=0.1$, one can note that the upper bound
of $|D|$ is rather insensitive to the choice of $d_0$, whereas the
lower bound in \eqref{eq:2.size-estim-EIT-less-conduct-fat-incl}
improves as $d_0$ increases. The converse situation occurs when
the inclusion is made by material of higher conductivity, see
Figures \ref{fig:T1_2D_pos}(b) and \ref{fig:T2_2D_pos}(b).

\begin{figure}[ht]
   \begin{minipage}{.48\textwidth}
     \centering
     \includegraphics[width=6cm]{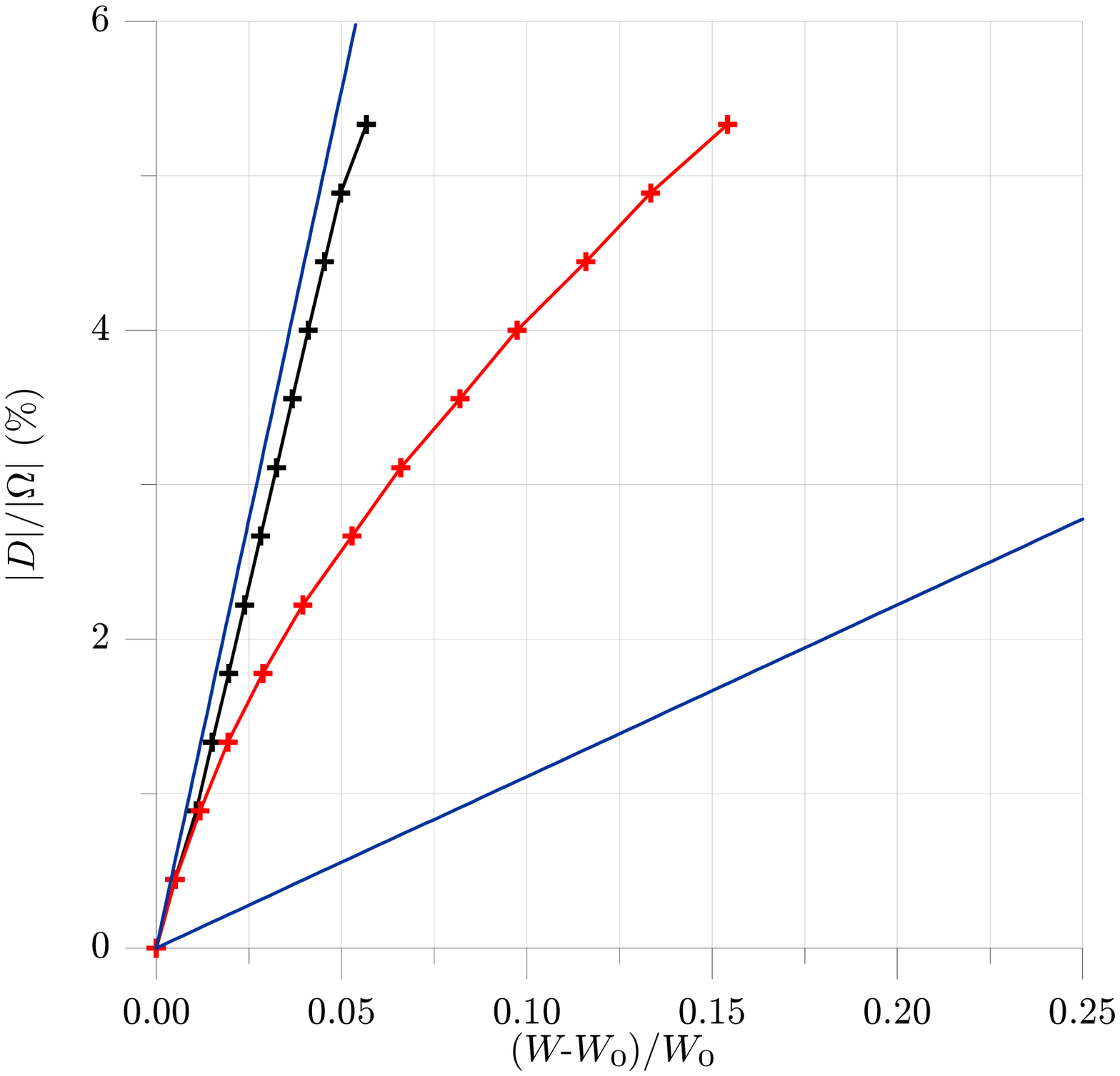}\\
     \centering{(a)}
   \end{minipage}
   \begin{minipage}{.48\textwidth}
     \centering
     \includegraphics[width=6cm]{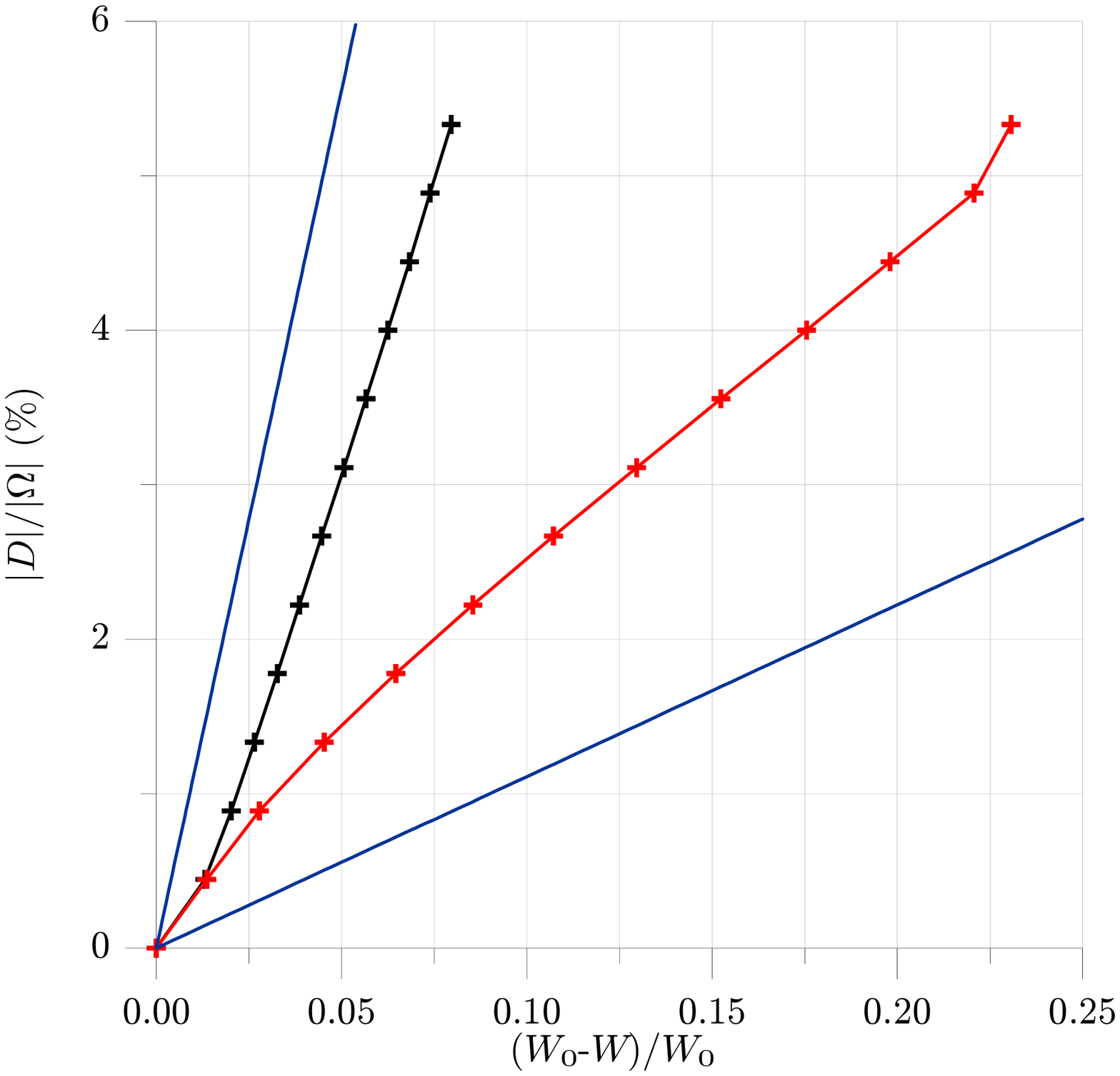}\\
     \centering{(b)}
   \end{minipage}
   \caption{Numerical size estimates for inclusions of general shape  
generated from a generic element inside $\Omega$ for test $T_1$ of  
Figure \ref{fig:Test numerici_2D}(a) ($21 \times 21$ FE mesh, $d_0=2 
$): $k=0.1$ (a), $k=10$ (b).}
   \label{fig:T1_2D_shp}
\end{figure}

\begin{figure}[ht]
   \begin{minipage}{.48\textwidth}
     \centering
     \includegraphics[width=6cm]{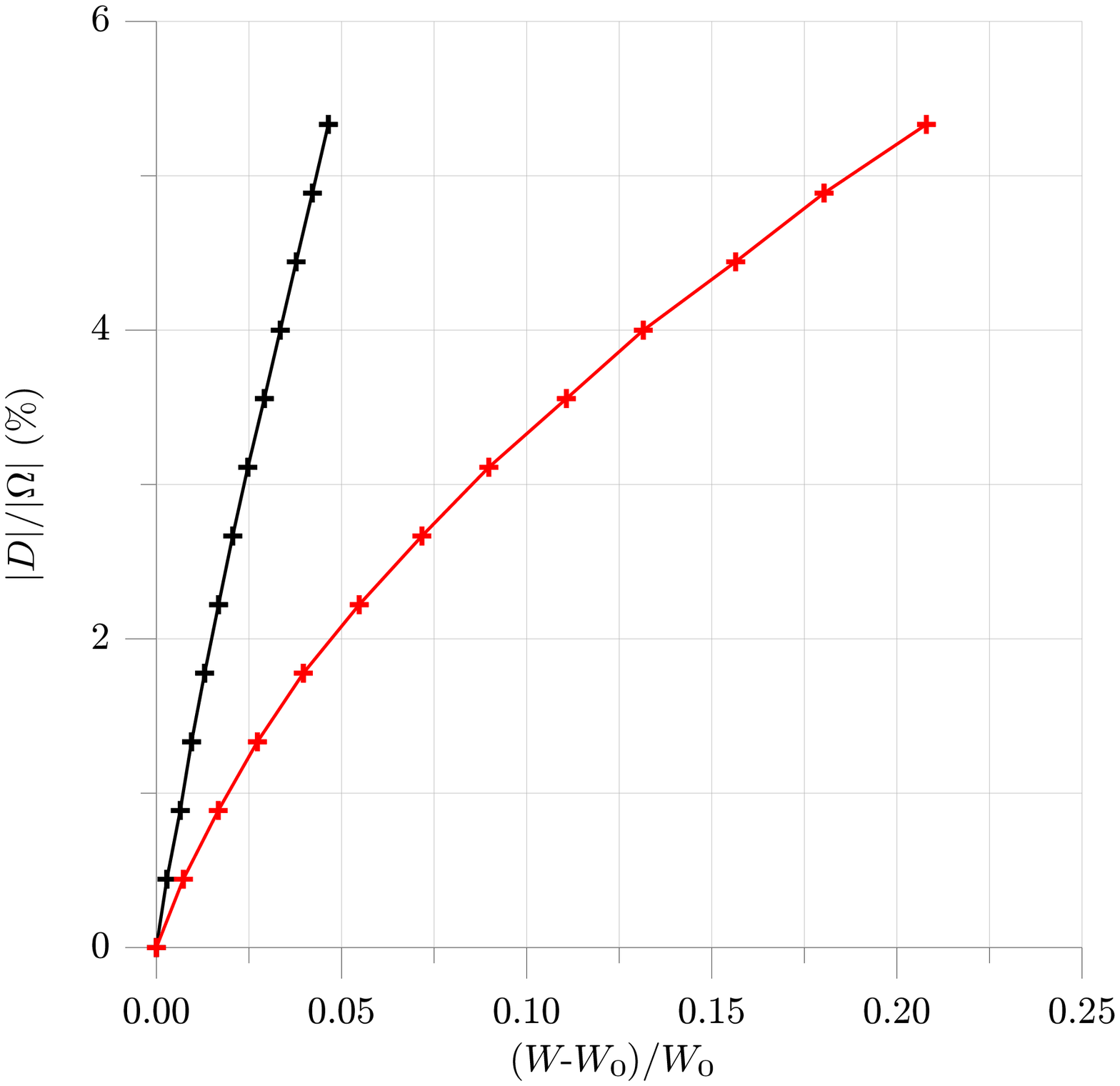}\\
     \centering{(a)}
   \end{minipage}
   \begin{minipage}{.48\textwidth}
     \centering
     \includegraphics[width=6cm]{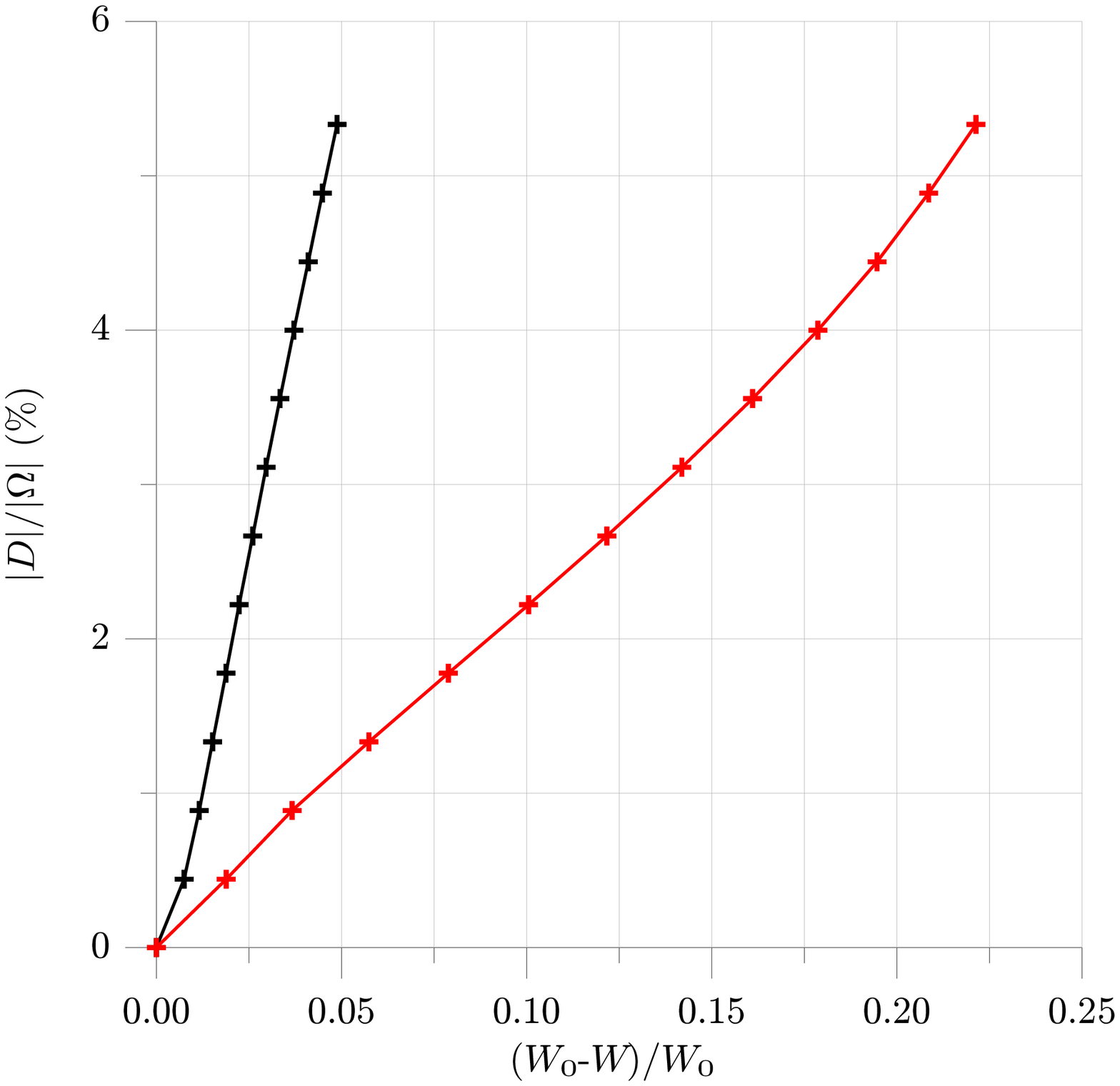}\\
     \centering{(b)}
   \end{minipage}
   \caption{Numerical size estimates for inclusions of general shape  
generated from a generic element inside $\Omega$ for test $T_2$ of  
Figure \ref{fig:Test numerici_2D}(b) ($21 \times 21$ FE mesh, $d_0=2 
$): $k=0.1$ (a), $k=10$ (b).}
   \label{fig:T2_2D_shp}
\end{figure}

As a second class of experiments, we have considered inclusions of
general shape on a FE mesh of $15 \times 15$ HC elements. More
precisely, each inclusion is the union of elements having at least
a common side and being at least $d_0=2$ elements far from the
boundary $\partial \Omega$. Results are collected in Figures
\ref{fig:T1_2D_shp} and \ref{fig:T2_2D_shp}.

The straight lines drawn in Figures \ref{fig:T1_2D_pos} and
\ref{fig:T1_2D_shp} correspond to the theoretical size estimates
for test $T_1$ of Figure \ref{fig:Test numerici_2D}(a). For both
cases $k=0.1$ and $k=10$ we have
\begin{equation}
     \label{eq:3.theor-size-T1}
   \frac{1}{9} \frac{|W-W_0|}{W_0} \leq \frac{|D|}{|\Omega|}\leq \frac 
{10}{9} \frac{|W-W_0|}{W_0}.
\end{equation}
The comparison with the region of the plane
$\left(\frac{|D|}{|\Omega|}, \frac{|W-W_0|}{W_0} \right)$ covered
by the corresponding numerical experiments confirms, as already
remarked in \cite{l:abfmrt04} in the context of linear elasticity,
that practical applications of the size estimates approach lead to
less pessimistic results with respect to those obtained via the
theoretical analysis.

\subsection{Three--dimensional case} \label{subsec:3D}

The first part of this subsection is devoted to the extension to
the 3--D case of the numerical simulations given in
\ref{subsec:2D}. In the second part, we shall investigate on the
effect of the oscillation character of the Neumann data on the
upper bound of size inclusion.

\begin{figure}[h]
   \begin{minipage}{.48\textwidth}
     \centering
     \includegraphics[width=5cm]{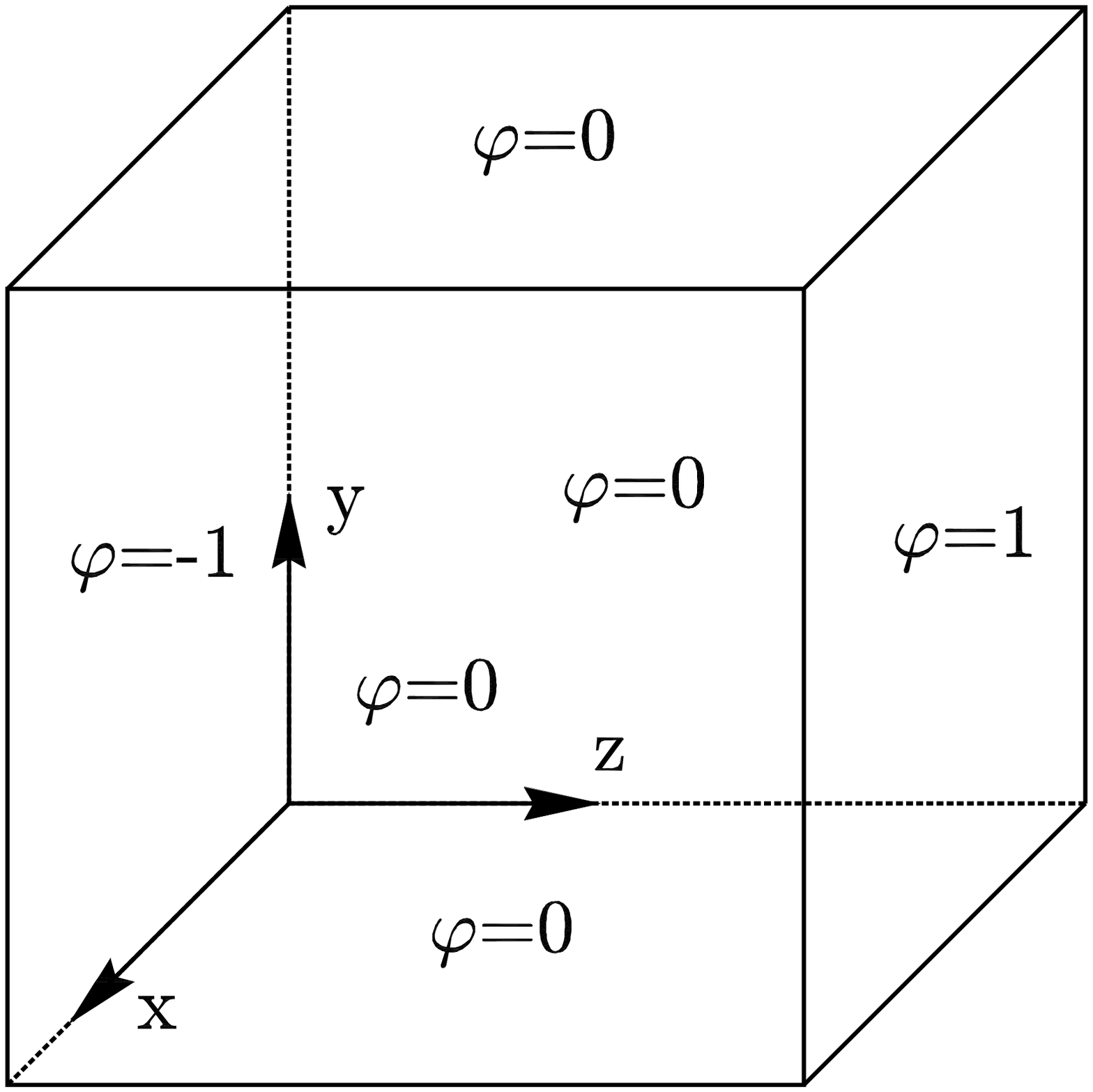}\\
     \centering{(a)}
   \end{minipage}
   \begin{minipage}{.48\textwidth}
     \centering
   \includegraphics[width=5cm]{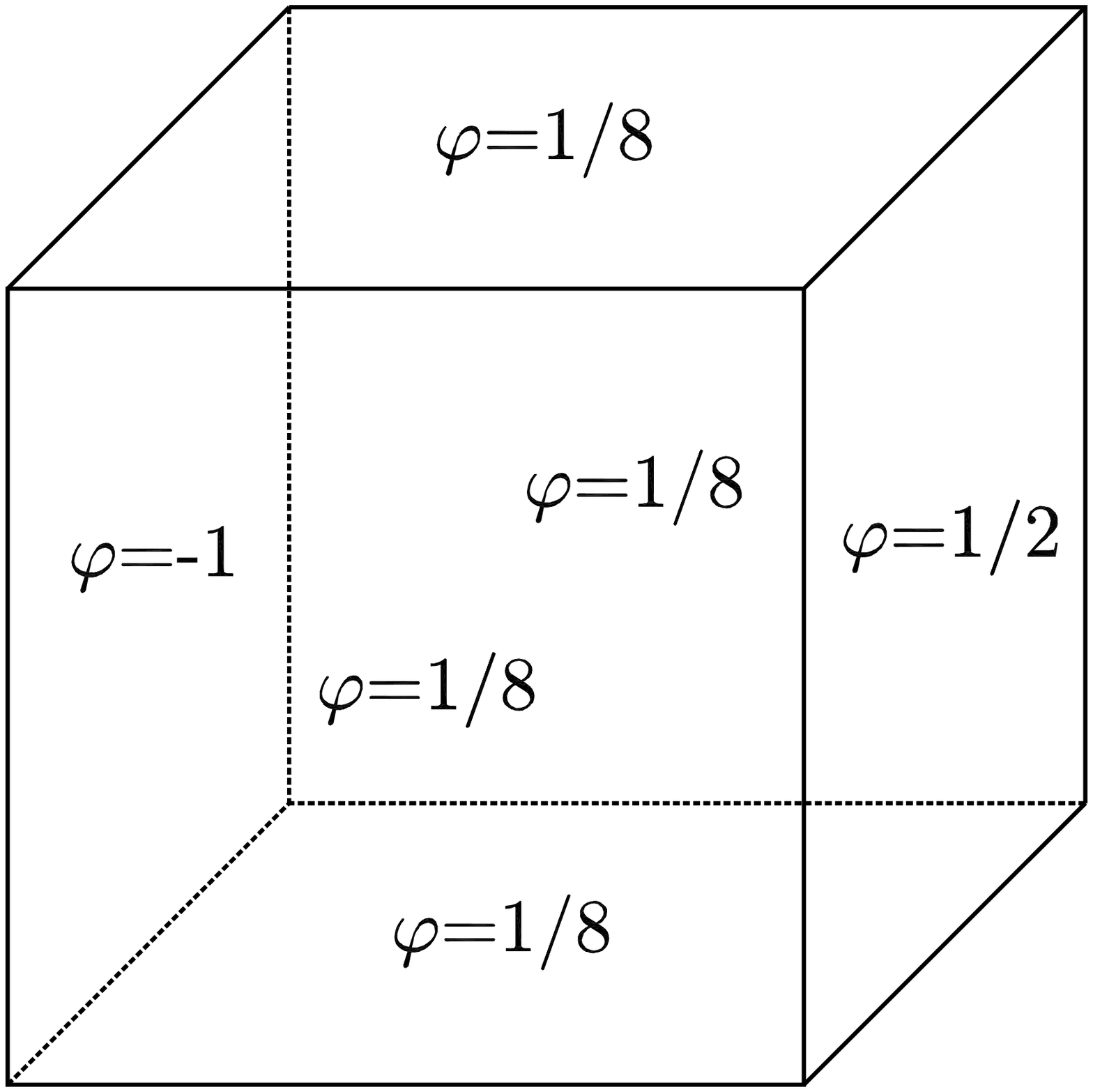}\\
     \centering{(b)}
   \end{minipage}
   \caption{Cubic conductor considered in 3--D numerical simulations  for the EIT model and applied current density fields: Test $T_1$ (a)  and Test $T_2$ (b).}
   \label{fig:Test numerici_3D}
\end{figure}

\begin{figure}[t]
   \begin{minipage}{.48\textwidth}
     \centering
     \includegraphics[width=6cm]{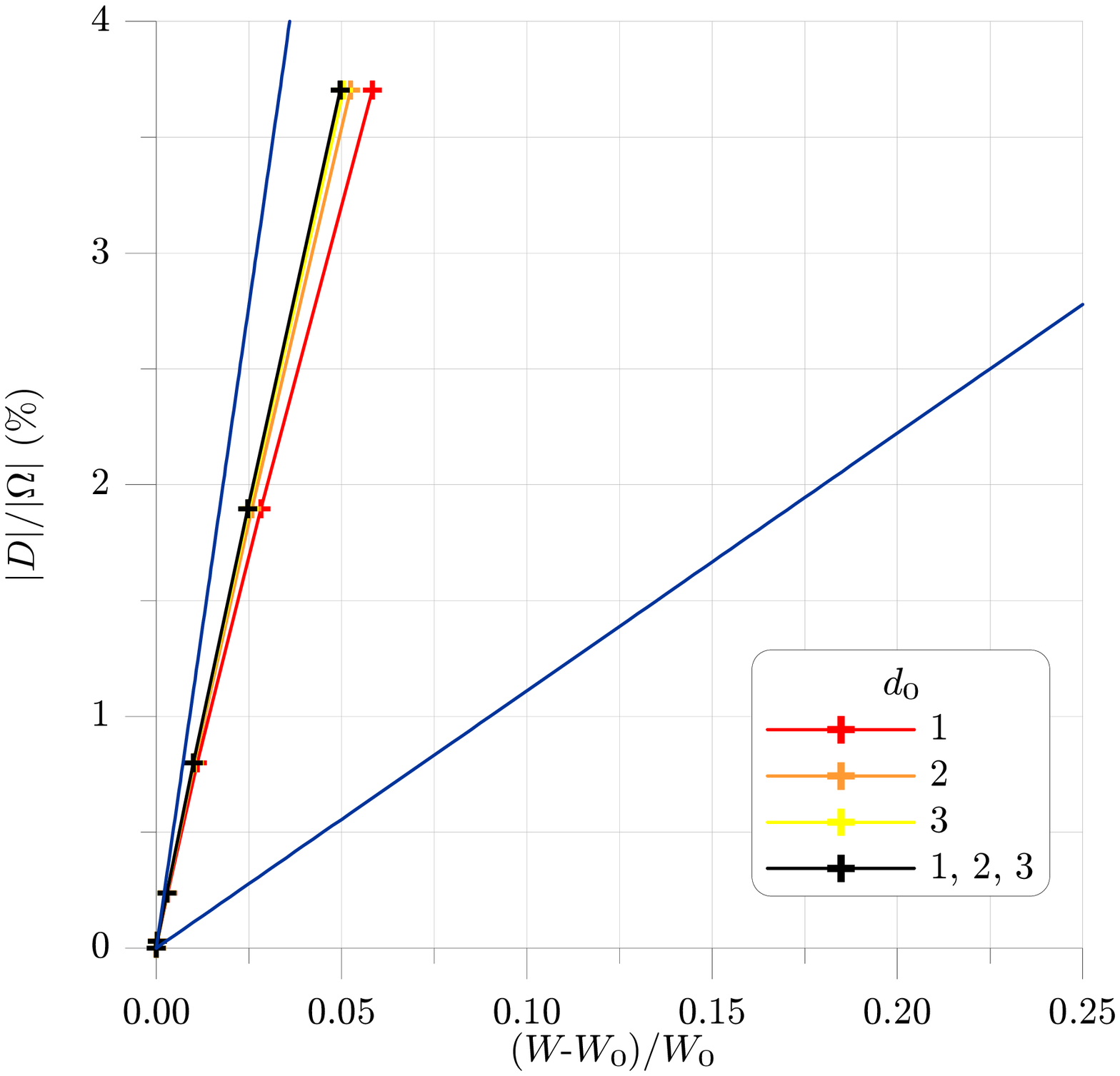}\\
     \centering{(a)}
   \end{minipage}
   \begin{minipage}{.48\textwidth}
     \centering
     \includegraphics[width=6cm]{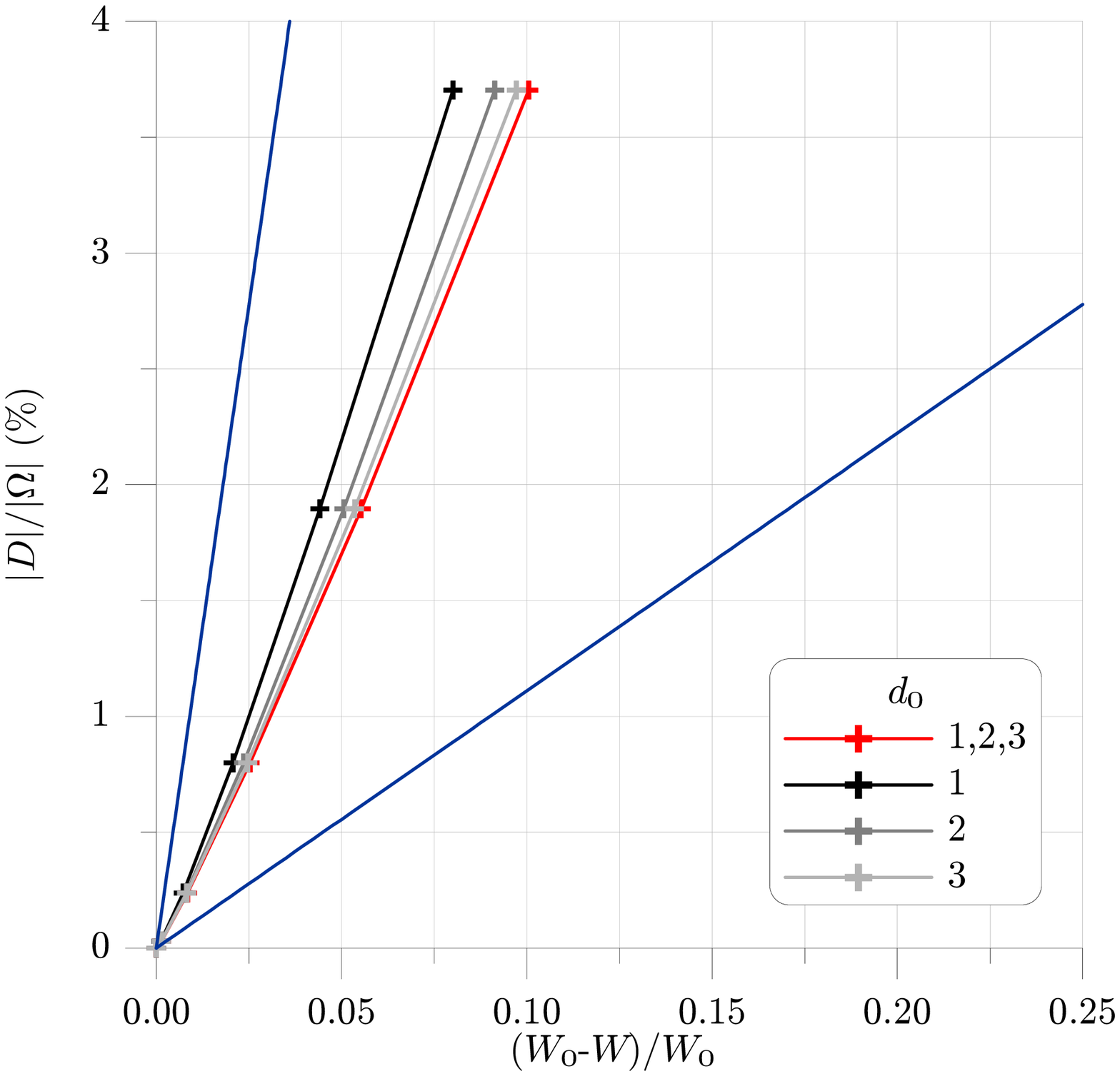}\\
     \centering{(b)}
   \end{minipage}
   \caption{Influence of $d_0$ for cubic inclusions in test $T_1$ of
Figure \ref{fig:Test numerici_3D} ($20 \times 20 \times 20$ FE
mesh): $k=0.1$ (a), $k=10$ (b).}
   \label{fig:T1_3D_pos}
\end{figure}

\begin{figure}[h]
   \begin{minipage}{.48\textwidth}
     \centering
     \includegraphics[width=6cm]{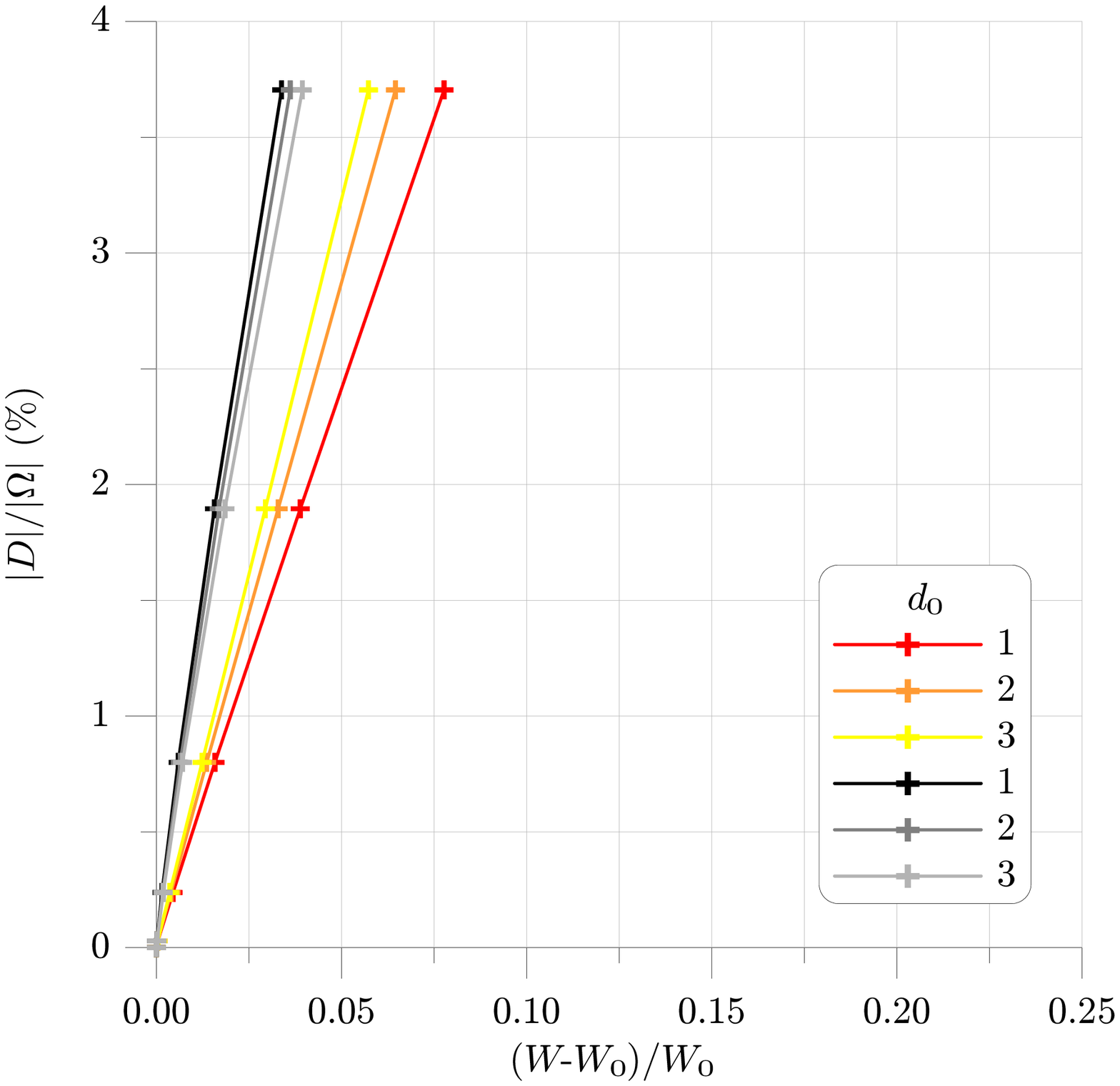}\\
     \centering{(a)}
   \end{minipage}
   \begin{minipage}{.48\textwidth}
     \centering
     \includegraphics[width=6cm]{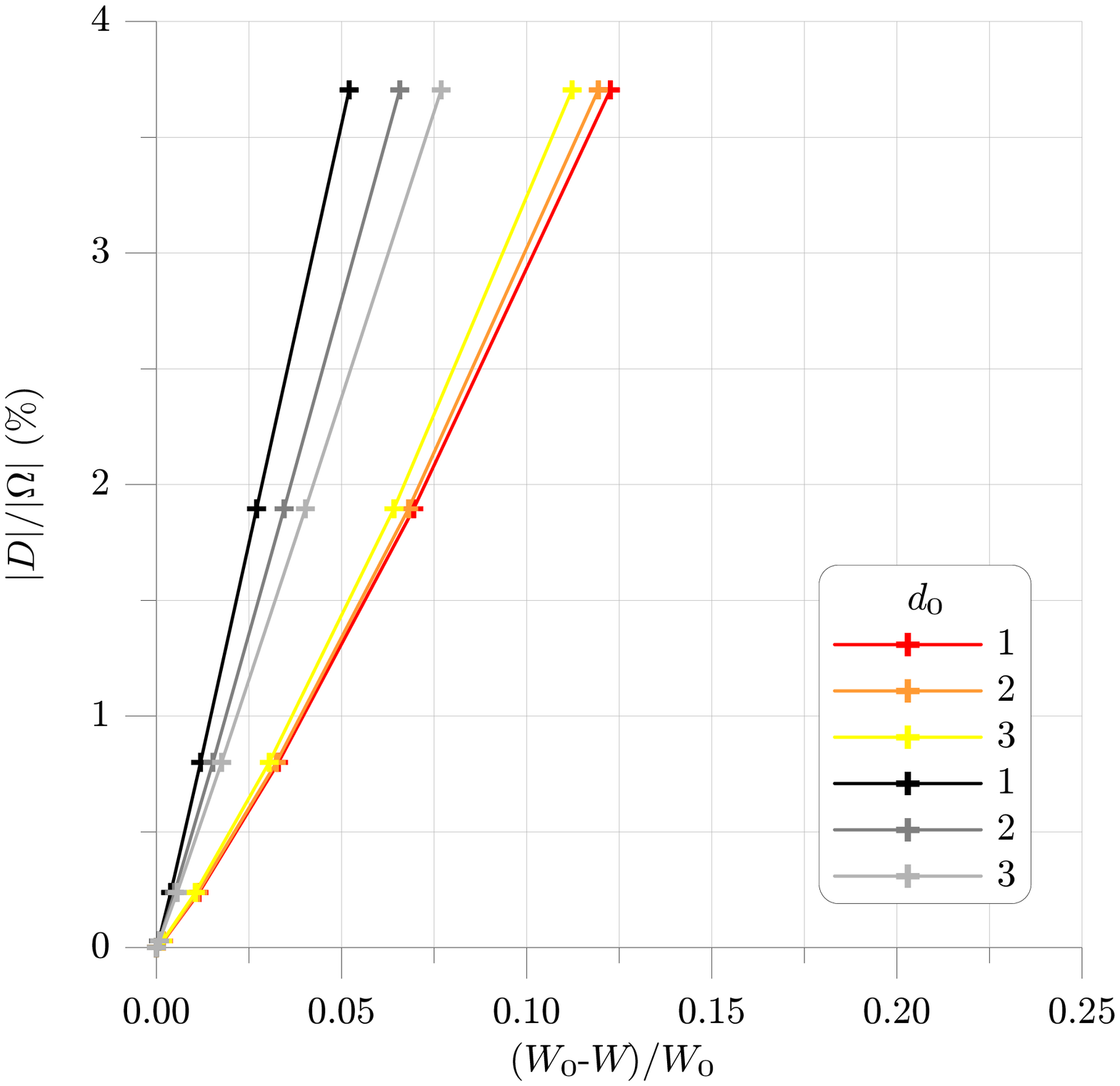}\\
     \centering{(b)}
   \end{minipage}
   \caption{Influence of $d_0$ for cubic inclusions in test $T_2$ of  
Figure \ref{fig:Test numerici_3D} ($20 \times 20 \times 20$ FE mesh):
$k=0.1$ (a), $k=10$ (b).}
   \label{fig:T2_3D_pos}
\end{figure}

Similarly to the 2--D case, a first series of numerical
simulations has been performed on an electrical conductor of cubic
shape, of side $l$, with the two current density fields
illustrated in Figure \ref{fig:Test numerici_3D}. In both cases, a
mesh of $20 \times 20 \times 20$ finite elements has been
considered when performing simulations in presence of cubic
inclusions. The results are illustrated in Figures
\ref{fig:T1_3D_pos} and \ref{fig:T2_3D_pos}. Figure
\ref{fig:T1_3D_pos} contains also the straight lines corresponding
to the theoretical size estimates for test $T_1$ of Figure
\ref{fig:Test numerici_3D}, that is
\begin{equation}
     \label{eq:3.theor-size-T1-3D}
   \frac{1}{9} \frac{|W-W_0|}{W_0} \leq \frac{|D|}{|\Omega|}\leq \frac 
{10}{9} \frac{|W-W_0|}{W_0}.
\end{equation}


In order to deal with inclusions of general shape, however, the
numerical experiments require some restrictions to reduce the
computer time. A rough estimate of the computational cost can be
obtained noting that the numerical effort is essentially due to
the decomposition of the matrix associated to the linear system
\eqref{eq:EIT_min_w} and to the computation of its solution.
Denoting by $m$ the number of the equations and by $b$ the half
bandwidth of the matrix, the decomposition requires $m(b-1)$
multiplications and $m b (b-1)$ additions, whereas the computation
of the solution involves $m b$ multiplications.

Therefore, for each given inclusion in a $20\times20\times20$ FE
mesh, a linear system of $10648$ ($b=1015$) equations has to be
solved, requiring a computer time of approximately $86$ s working
on an Opteron $2.4$ GHz computer. Since the number of all possible
inclusions formed by $n_i$ elements on a mesh of $n_e \times n_e
\times n_e$ is $\frac{n_e^3 !}{n_i ! (n_e^3-n_i)!}$, the way to
calculate all the possible case is practically impossible. Indeed
by considering that the $20\times20\times20$ is formed by $8000$
elements and that, if the ratio $|D|/|\body|$ is less than $6\%$
that is $480$ elements, the number of cases to analyze is
$69.1183\times 10^{785}$.

\begin{figure}[ht]
   \begin{minipage}{.48\textwidth}
     \centering
     \includegraphics[width=6cm]{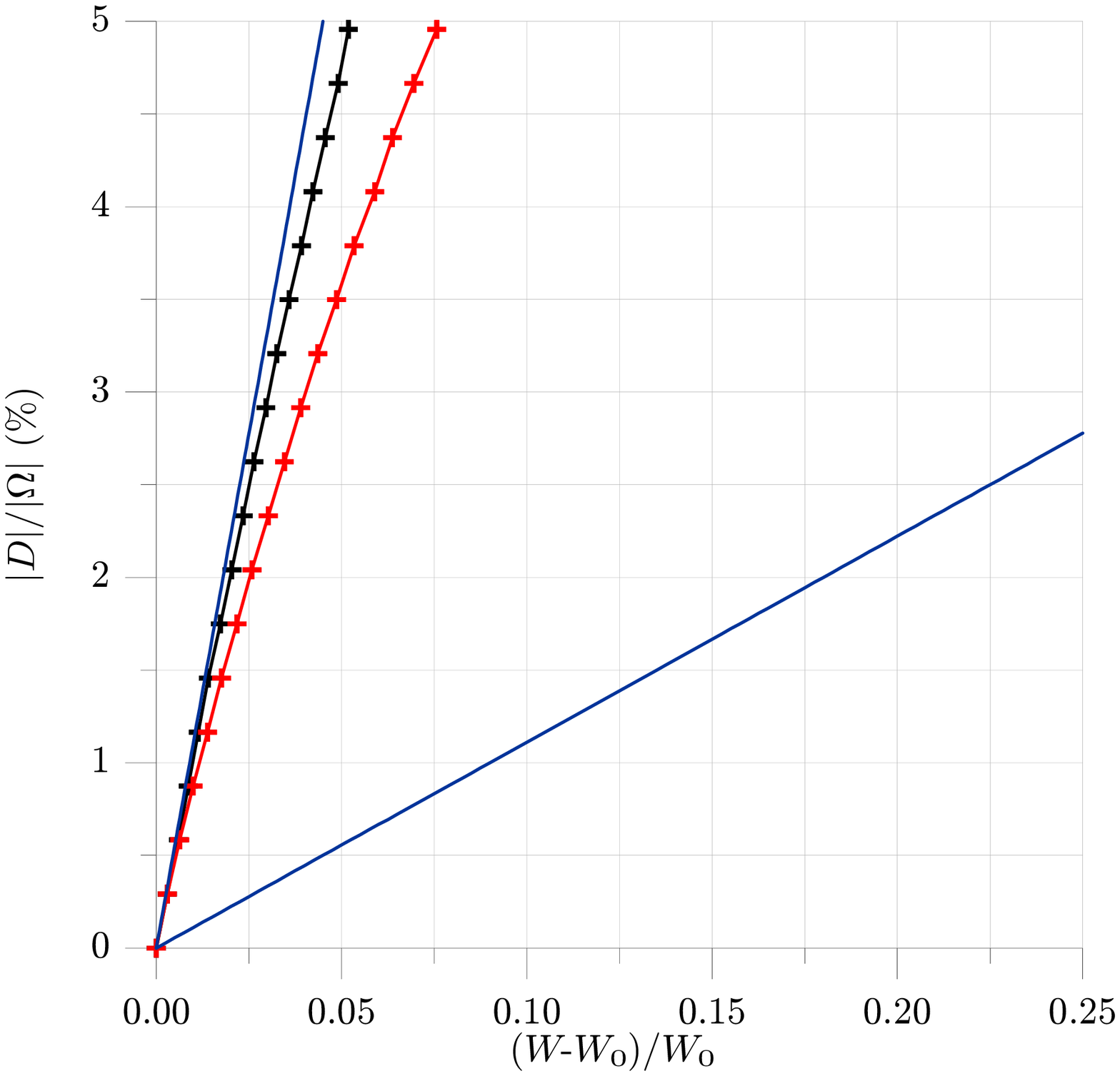}\\
     \centering{(a)}
   \end{minipage}
   \begin{minipage}{.48\textwidth}
     \centering
     \includegraphics[width=6cm]{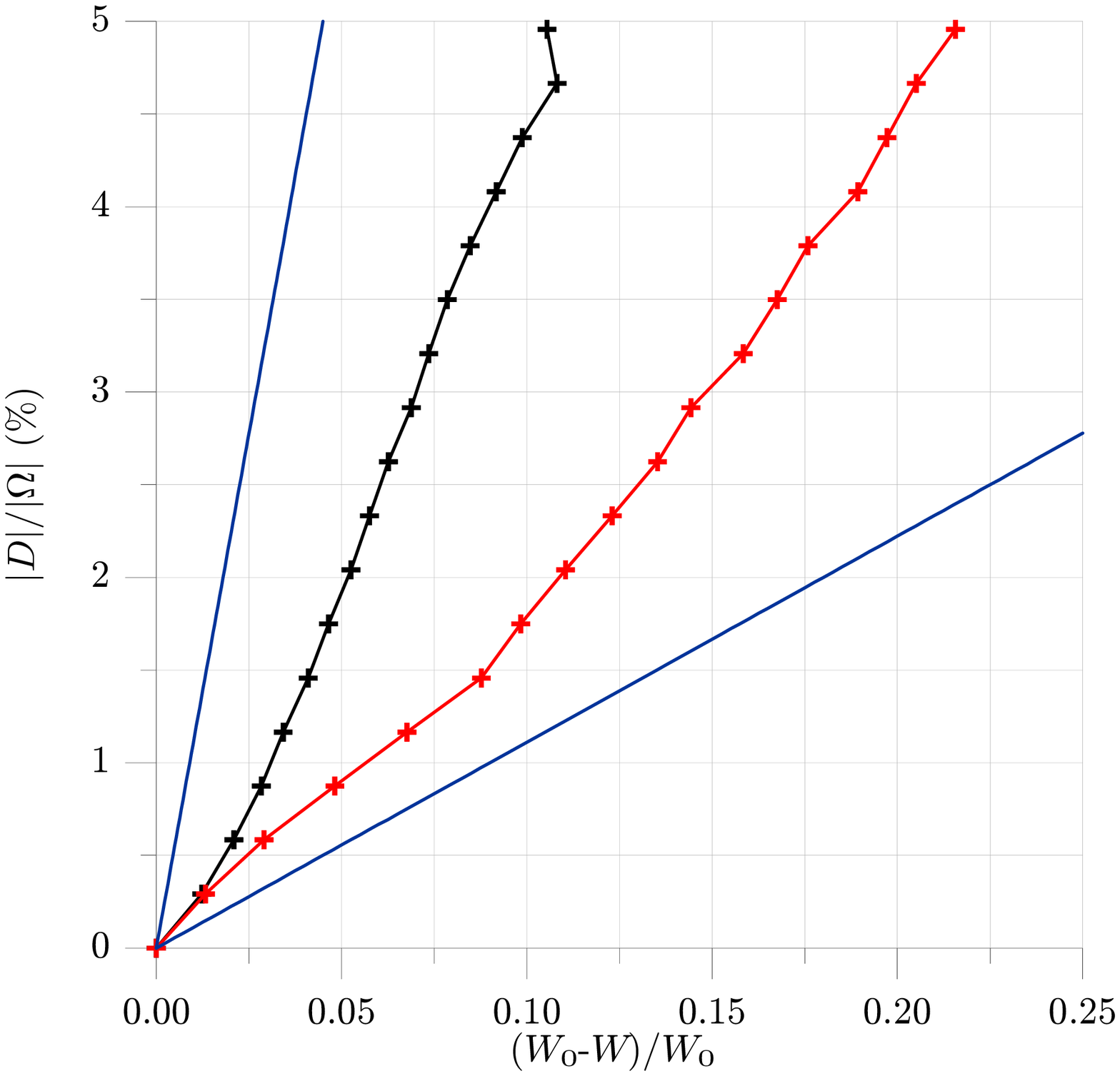}\\
     \centering{(b)}
   \end{minipage}
   \caption{Numerical size estimates for inclusions of general shape  
generated from a generic element belonging to an eight of the cube  
for test $T_1$ of Figure \ref{fig:Test numerici_3D}(a) ($7 \times 7  
\times 7$ FE mesh, $d_0=1 $): $k=0.1$ (a), $k=10$ (b).}
   \label{fig:T1_3D_shp}
\end{figure}

\begin{figure}[ht]
   \begin{minipage}{.48\textwidth}
     \centering
     \includegraphics[width=6cm]{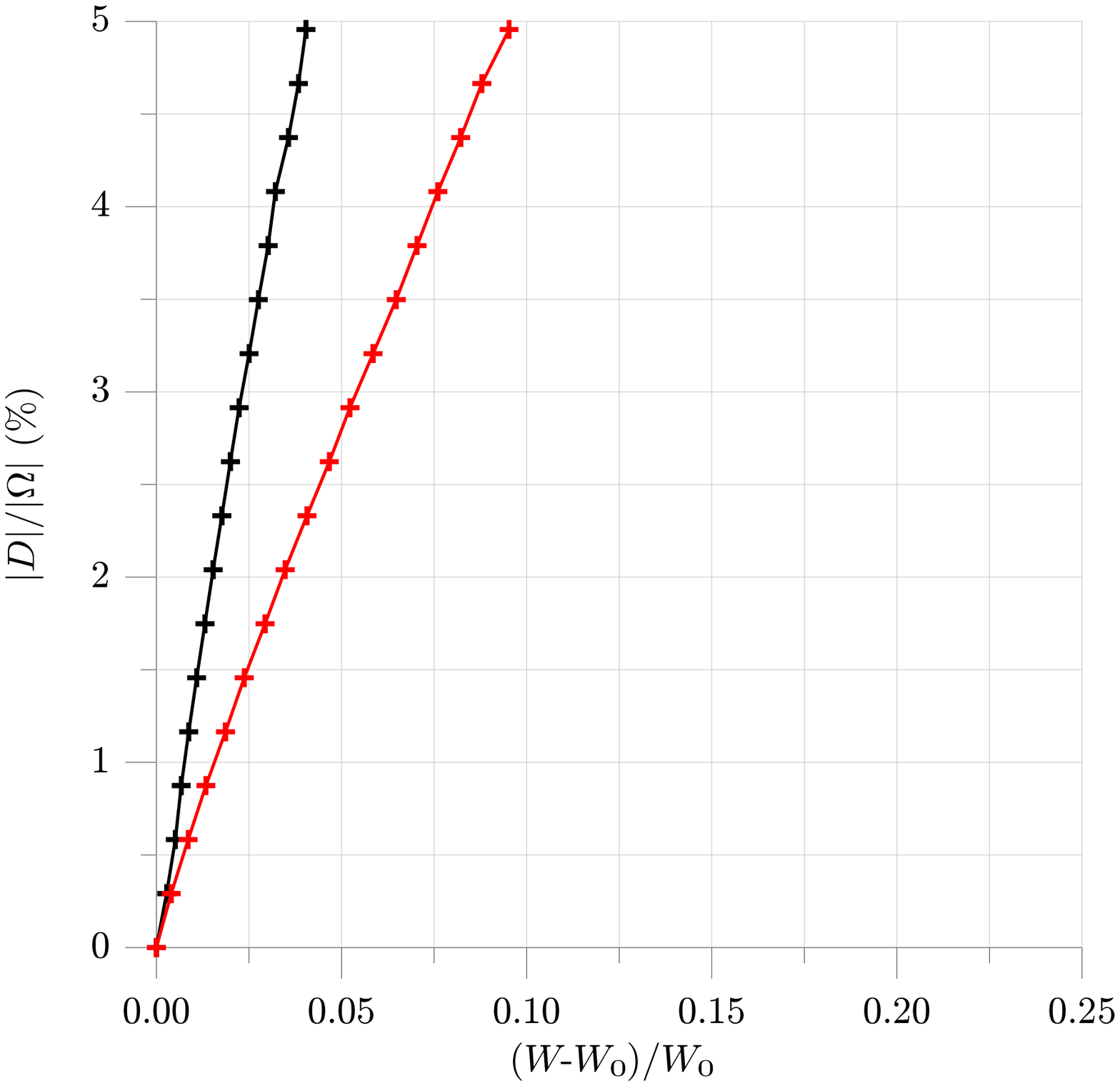}\\
     \centering{(a)}
   \end{minipage}
   \begin{minipage}{.48\textwidth}
     \centering
     \includegraphics[width=6cm]{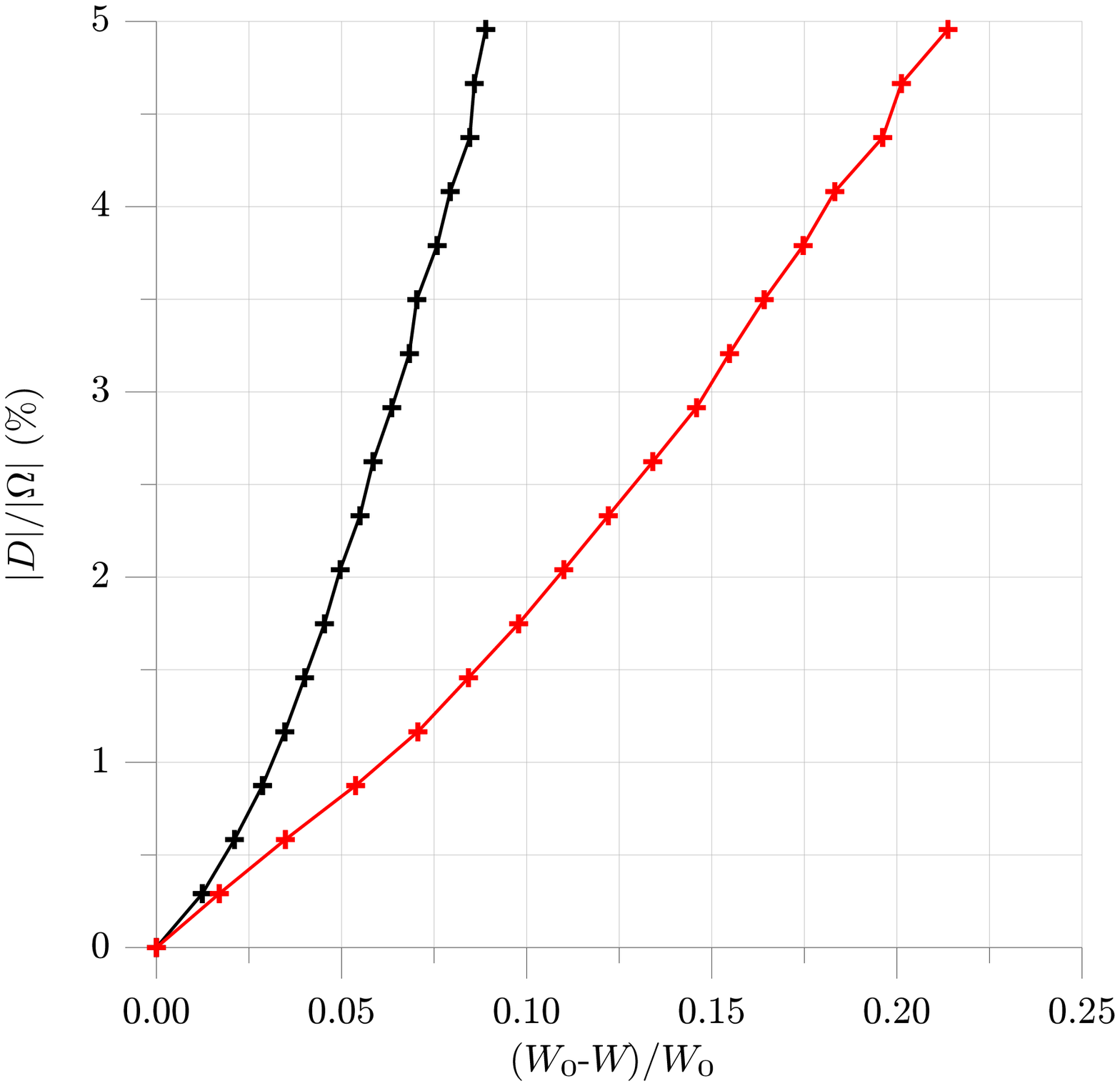}\\
     \centering{(b)}
   \end{minipage}
   \caption{Numerical size estimates for inclusions of general shape  
generated from a generic element belonging to an eight of the cube  
for test $T_2$ of Figure \ref{fig:Test numerici_3D}(a) ($7 \times 7  
\times 7$ FE mesh, $d_0=1 $): $k=0.1$ (a), $k=10$ (b).}
   \label{fig:T2_3D_shp}
\end{figure}

In order to reduce the computer time significantly we have
considered a $7\times 7 \times 7$ mesh generating a system of
$729$ equations. Despite of this, the number of possible cases to
consider still remains very high; for instance, for inclusions
formed by $5$ elements, one should solve about $3.8\times 10^{10}$
linear systems. Therefore, we decided to restrict our analysis to
inclusions satisfying the following additional hypotheses:
\begin{enumerate}[i)]
   \item the inclusion is the union of elements having at least one
    common face and it is formed by starting from a generic element  
inside an
    octant of the cube (this last assumption is not really  
restrictive due to the symmetries of the problem);
   \item $d_0=1$.
\end{enumerate}
For inclusions formed by $1,...,7$ elements, we have considered
all possible inclusions satisfying the limitations $i)$ and $ii)$,
whereas for inclusions formed by $8,...,17$ elements we have
considered a random sample because of the high computational cost.
For these cases, the ratio between the sample dimension and that
of all the data approximately spans between $20\%$ for inclusions
formed by $8$ elements and $0.01\%$ for inclusions formed by $17$
elements. The results are presented in Figures \ref{fig:T1_3D_shp}
and \ref{fig:T2_3D_shp} for Test $T_1$ and Test $T_2$,
respectively. In Figure \ref{fig:T1_3D_shp}, the straight lines
corresponding to the theoretical bounds
\eqref{eq:3.theor-size-T1-3D} for Test $T_1$ are also drawn. As
already remarked in the treatment of the 2--D case, the
theoretical analysis leads to rather pessimistic results with
respect to those obtained by the numerical simulations, especially
when the inclusion is softer than the surrounding material.


The Neumann data considered in the above experiments give raise to
potential fields inside the conductor with nonvanishing gradient.
In the general case, when the gradient of the solution may vanish,
we expect, accordingly to Theorems
\ref{theo:size-estim-EIT-general},
\ref{theo:size-estim-EIT-fat-incl}, that the upper bounds
deteriorate as the frequency $F[\varphi]$ given by
\eqref{eq:2.frequency} increases. Since $F[\varphi]$ is a ratio
which measures the frequency of oscillation of $\varphi$, we are
interested to investigate on the effectiveness of size estimates
approach for oscillating Neumann data.

In particular, the numerical simulations have been carried out for
the cubic electrical conductor considered before and choosing the
following Neumann data:
\begin{equation} \label{eq:3.oscill-Neu-data}
     \left.
         \begin{array}{crl}
             \varphi = & - \cos \frac{n \pi x}{l}  & \qquad \hbox 
{on } z=0 , \\
             \varphi = &  \cos \frac{n \pi x}{l}   & \qquad \hbox 
{on } z=l , \\
             \varphi = & 0                                   & \qquad  
\hbox{elsewhere on } \partial \Omega,
         \end{array}
     \right\}
     \quad \hbox{for } n=0, 1, 2 .
\end{equation}

Case $n=0$ has been already discussed at the beginning of this
paragraph and corresponds to the simple case in which the gradient
of the unperturbed solution $u_0$ does not vanish in $\Omega$.

The two other cases are examples of Neumann data with higher
frequency $F[\varphi]$. More precisely, the corresponding
solutions $u_0$ have critical lines of equation 
\[ 
\left\{ x= \frac{l}{n} \left( \frac{1}{2}+i \right),  z= \frac{l}{n} \left( \frac{1}{2}+j \right) \right\}, \quad i,j = 0, ..., n-1.
\]

\begin{figure}[h]
\begin{minipage}{0.49\textwidth}
   \includegraphics[height=6cm]{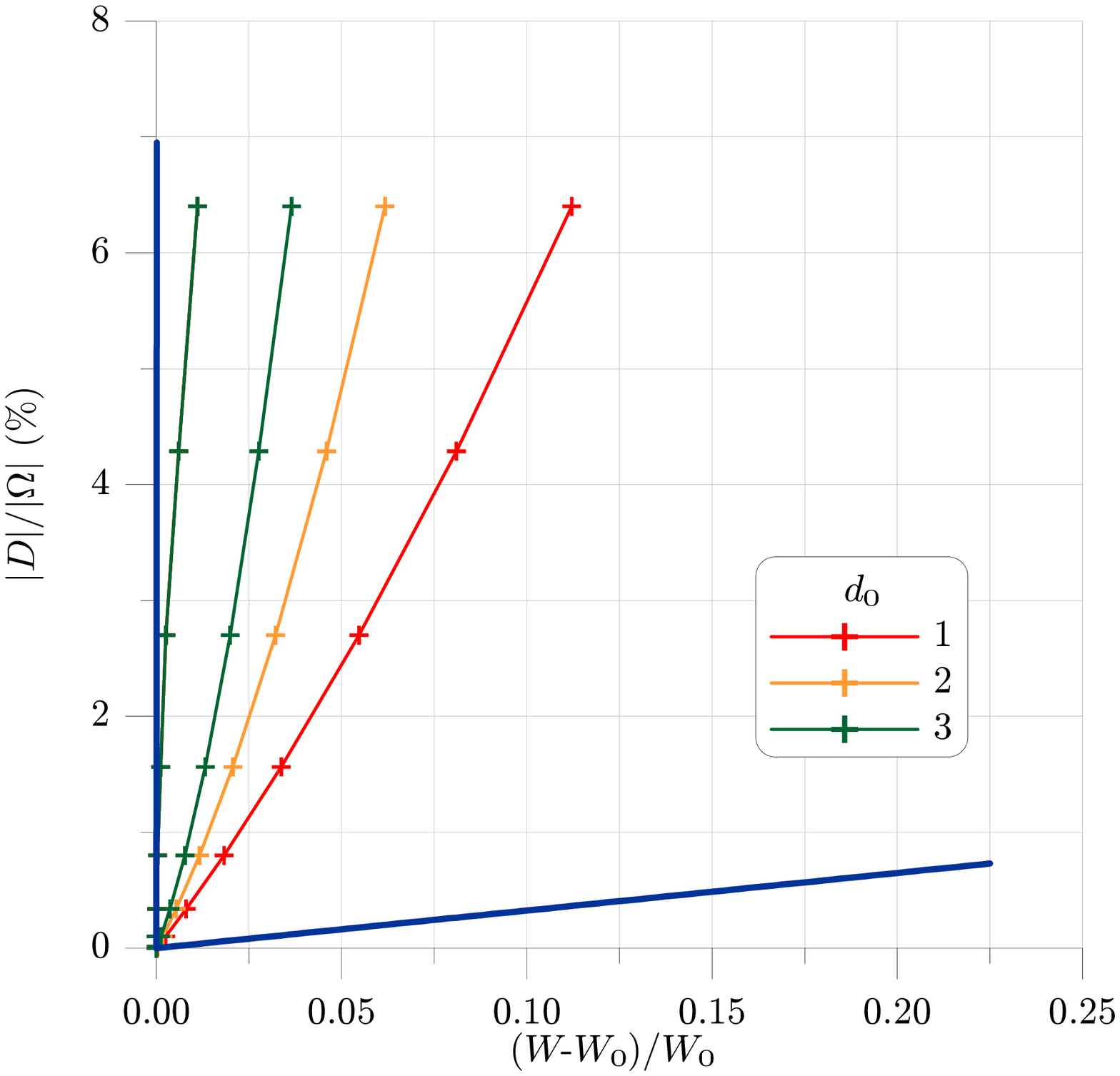} \\
   \centering{(a)}
\end{minipage}
\begin{minipage}{0.49\textwidth}
   \includegraphics[height=6cm]{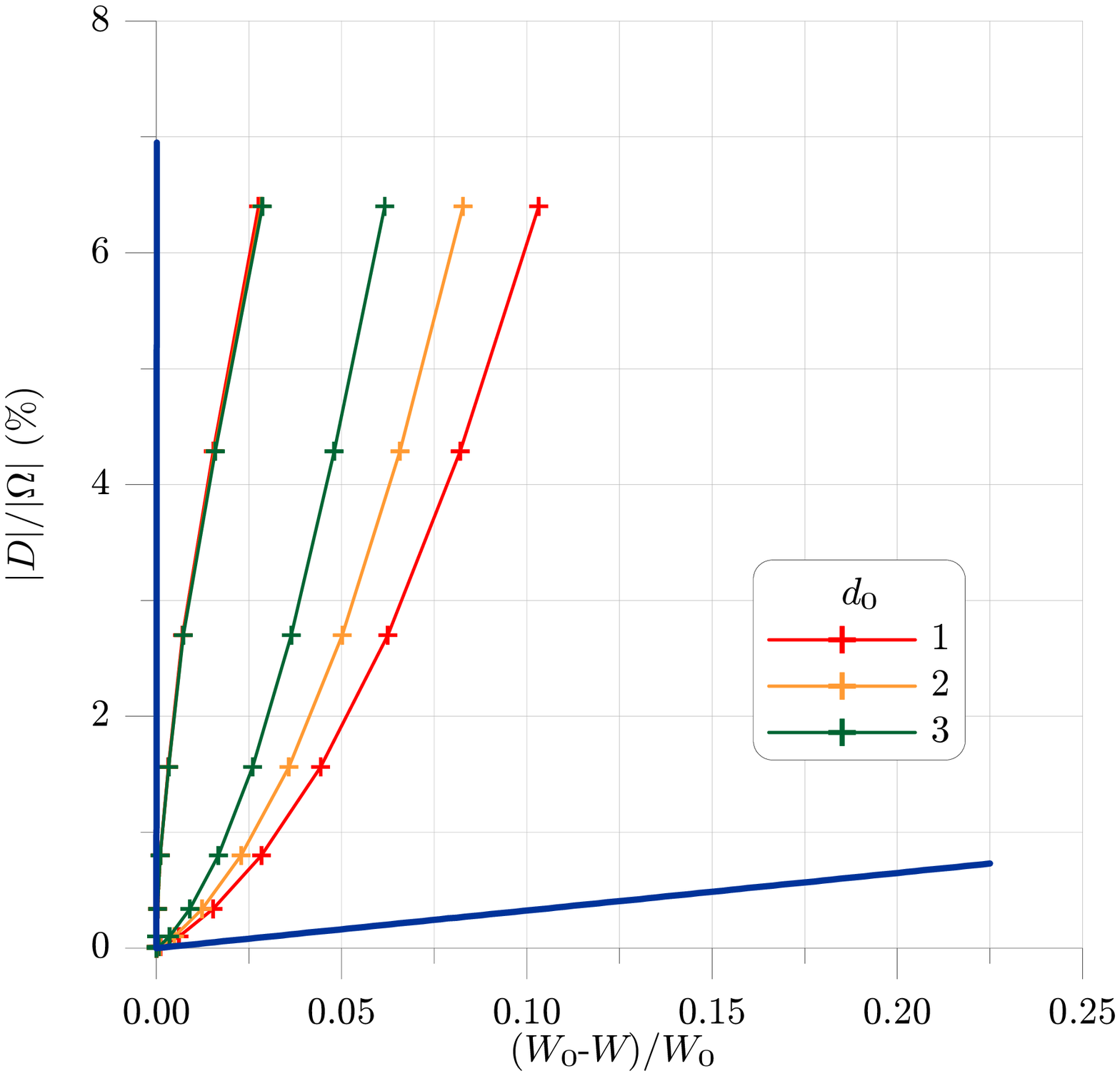} \\
   \centering{(b)}
\end{minipage}
\caption{Cubic electrical conductor with Neumann data as in case
$n=1$ of \eqref{eq:3.oscill-Neu-data}: lower and upper bound of
the power gap for different values of $d_0$ ($k=0.1$ (a) and
$k=10$ (b)) on a $20 \times 20 \times 20$ mesh.} \label{fig:cos_1}
\end{figure}

\begin{figure}[h]  
\begin{minipage}{0.49\textwidth}
   \includegraphics[height=6cm]{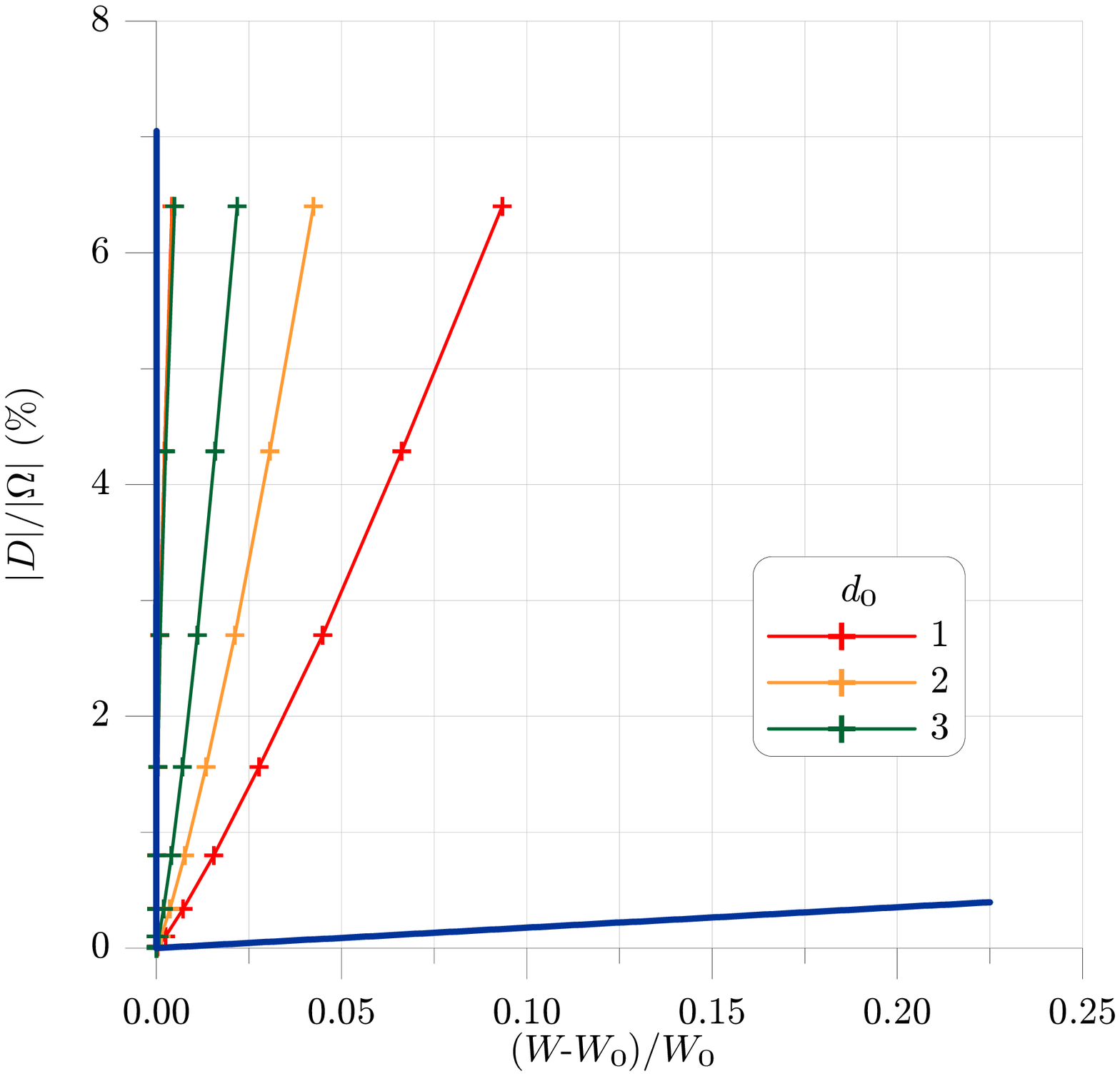} \\
   \centering{(a)}
\end{minipage}
\begin{minipage}{0.49\textwidth}
   \includegraphics[height=6cm]{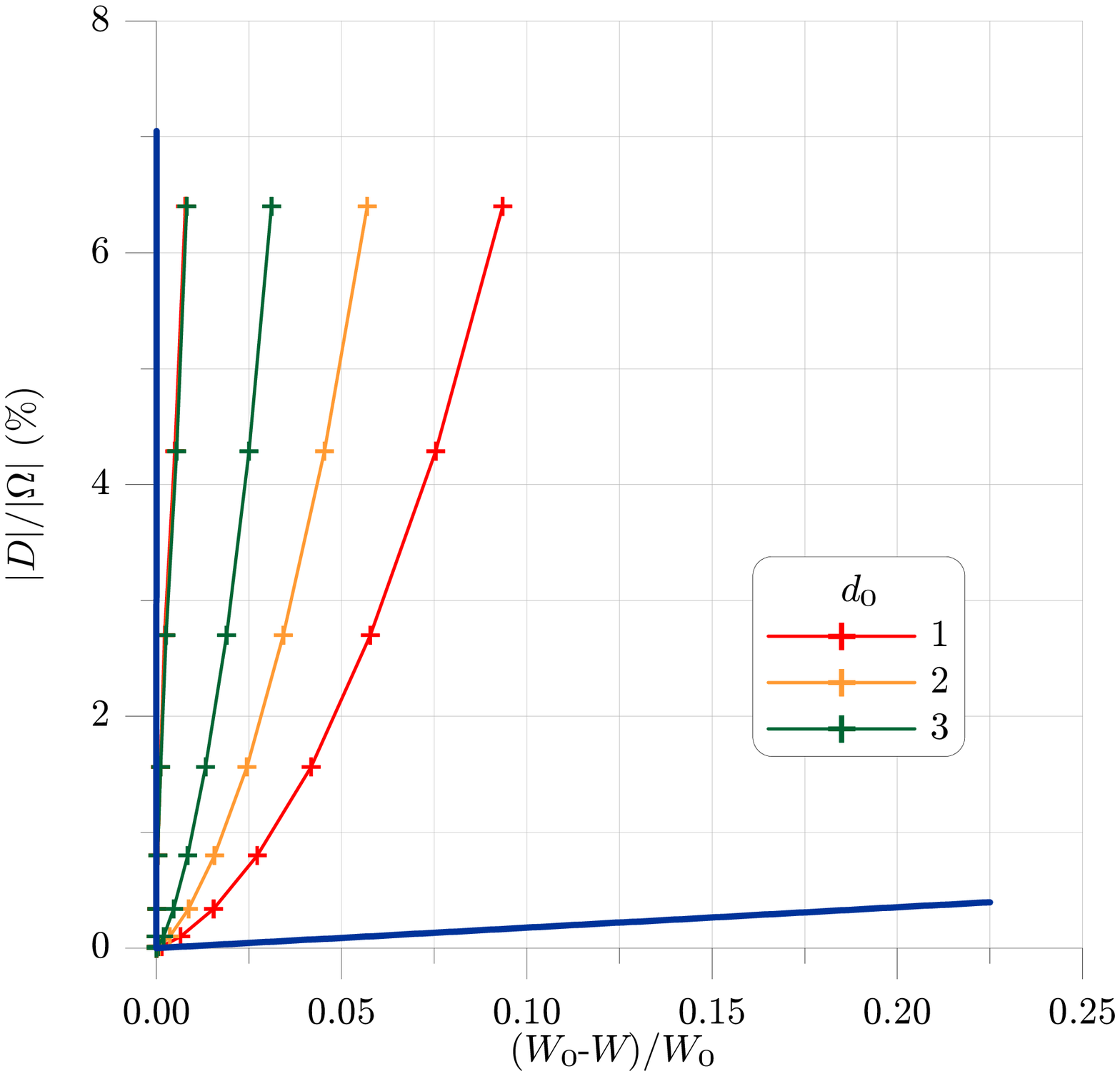} \\
   \centering{(b)}
\end{minipage}
\caption{Cubic electrical conductor with Neumann data as in case
$n=2$ of \eqref{eq:3.oscill-Neu-data}: lower and upper bound of
the power gap for different values of $d_0$ ($k=0.1$ (a) and
$k=10$ (b)) on a $20 \times 20 \times 20$ mesh.} \label{fig:cos_2}
\end{figure}

The mesh employed is made by $20 \times 20 \times 20$ HC finite
elements. The analysis has been focussed on cubic inclusions
having volume up to $6 \%$ of the total volume of the specimen and
conductivity $k=0.1$ and $k=10$. The numerical results in case
$n=1$ and $n=2$ are presented in Figures \ref{fig:cos_1} and
\ref{fig:cos_2}, respectively. The numerical results show that the
lower bound in size estimates
\eqref{eq:2.size-estim-EIT-more-conduct-fat-incl},
\eqref{eq:2.size-estim-EIT-less-conduct-fat-incl} improves as
$d_0$ increases, whereas the upper bound of $|D|$ is rather
insensitive to the choice of $d_0$.

Theoretical estimates for cases $n=1$ and $n=2$ of
\eqref{eq:3.oscill-Neu-data} are given by
\begin{equation}\label{eq:cos_est}
     \begin{split}
       \hbox{for } & k>1:  \\
        & \frac{\tanh \frac{n\pi}{2}}{ n\pi(k-1) } \frac{W_0-W}{W_0}
       \leq
       \frac{|D|}{|\Omega|}
       \leq
       \frac{1}{C_n} \frac{k}{k-1} \frac{\tanh \frac{n\pi}{2}}{n\pi}  
\frac{W_0-W}{W_0};  \\
       \hbox{for } & k<1:  \\
       &\frac{k}{n\pi(1-k)} \tanh \frac{n\pi}{2} \frac{W-W_0}{W_0}
       \leq
       \frac{|D|}{|\Omega|}
       \leq \frac{1}{C_n} \frac{1}{1-k} \frac{\tanh \frac{n\pi}{2}}{n 
\pi}
       \frac{W-W_0}{W_0},
     \end{split}
\end{equation}
where
$$
C_n= \frac{10}{n\pi \cosh^2 \frac{n\pi}{2} } \left ( \sinh
\frac{n\pi}{20} - \sin \frac{n\pi}{20} \right ), \quad \quad
n=1,2.
$$
The theoretical estimates are indicated in Figures \ref{fig:cos_1}
and \ref{fig:cos_2}. The slope of the straight line corresponding
to the upper bound is so high that it practically coincides with
the vertical axis, at least for the portion of graph near the
origin considered in this study. The theoretical lower bound
gives, for a fixed power gap, values significantly less than those
obtained in the numerical experiments.

\section{Numerical simulations for the complete EIT model} \label{sec: num-phys-EIT}

\subsection{Numerical model}  \label{subsec:nummodel-phys}

In this case, by using the same notation introduced in Section
\ref{sec:num-EIT}, the energy functional
$J:H^1(\body)\times\mathbb{R}^L \rightarrow \mathbb{R}$ related to
the variational formulation of problem
(\ref{eq:2.Phys-Neumann_pbm_with_incl}) is given by
\begin{equation} \label{eq:Phis_EIT_var_form}
    J(u, U^l) =  \frac{1}{2} \int_\body (1 + (k-1) \chi_D) \nabla u  
\cdot \nabla u +
    \frac{1}{2} \sum_{l=1}^{L} \frac{1}{z_l} \int_{\boundary_l} (u- 
U^l)^2 - \sum_{l=1}^{L} I_l U^l .
\end{equation}

Using HC interpolation for the potential field $u$ and with the
notation introduced in Section \ref{sec:num-EIT}, the discrete
energy functional becomes
\begin{equation} \label{eq:Phis_EIT_var_form_discr}
   \begin{split}
    J(\mathbf{w}_e, U^l) & = \frac{1}{2} \sum_e \int_{\body_e} (1 +  
(k-1) \chi_D) (\nabla \mathbf{N}_e \mathbf{w}_e) \cdot   (\nabla  
\mathbf{N}_e \mathbf{w}_e) + \\
                         & + \frac{1}{2}  \sum_{l=1}^{L} \frac{1} 
{z_l} \sum_{\hat{e}} \int_{(\boundary_l)_e} (\mathbf{N}_{e} \mathbf{w} 
_{e} - U^l)^2  - \sum_{l=1}^{L} I_l U^l ,
   \end{split}
\end{equation}
or
\begin{equation} \label{eq:Phis_EIT_compact}
\begin{split}
    J(\mathbf{w}_e, U^l) & = \frac{1}{2} \sum_e \mathbf{w}_e^T \mathbf 
{K}_e \mathbf{w}_e + \\
                                         & + \frac{1}{2}  \sum_{l=1}^ 
{L} \frac{1}{z_l} \sum_{\hat{e}} ( \mathbf{w}_e^T \mathbf{K}_{ll}  
\mathbf{w}_e + (U^l)^2 - 2 \mathbf{w}_e^T \mathbf{K}_{el} U^l)  -  
\sum_{l=1}^{L} I_l
                                         U^l,
\end{split}
\end{equation}
having used  the compact notation
\begin{equation}\label{eq:Phis_EIT_mat}
     \begin{split}
       \mathbf{K}_e = & \int_{\body_e} (1 + (k-1) \chi_D) (\nabla  
\mathbf{N}_e)^T \nabla \mathbf{N}_e, \\
       \mathbf{K}_{ll} = & \int_{(\boundary_l)_e}  \mathbf{N}_e^T  
\mathbf{N}_e, \\
       \mathbf{K}_{el} = & \int_{(\boundary_l)_e}  \mathbf{N}_e^T.
     \end{split}
\end{equation}
We remark that the second sum in the right hand side of
\eqref{eq:Phis_EIT_var_form_discr} and
\eqref{eq:Phis_EIT_compact}, that  on $\hat{e}$, is extended only
to the elements under the electrodes.

Collecting the unknown parameters representing the potential field
in $\mathbf{w}$, those of the electrodes in $\mathbf{U}$ and the
current pattern in $\mathbf{I}$, by a standard method of
assembling we obtain the following linear system
\begin{equation}\label{eq:Phis_EIT_system_mat_form}
\begin{bmatrix}
   \mathbf{K}_{ww} & -\mathbf{K}_{wU} \\
   -\mathbf{K}_{wU}^T & \mathbf{K}_{UU} \\
\end{bmatrix} \begin{bmatrix}
                 \mathbf{w} \\
                 \mathbf{U} \\
               \end{bmatrix} = \begin{bmatrix}
                                 \mathbf{0} \\
                                 \mathbf{I} \\
                               \end{bmatrix} ,
\end{equation}
which can be efficiently solved taking advantage of the particular
structure of coefficient matrix.

\subsection{Results for 3--D cases} \label{subsec:phys-3D}

\begin{figure}[h]
   \begin{minipage}{.30\textwidth}
     \centering
     \includegraphics[width=4cm]{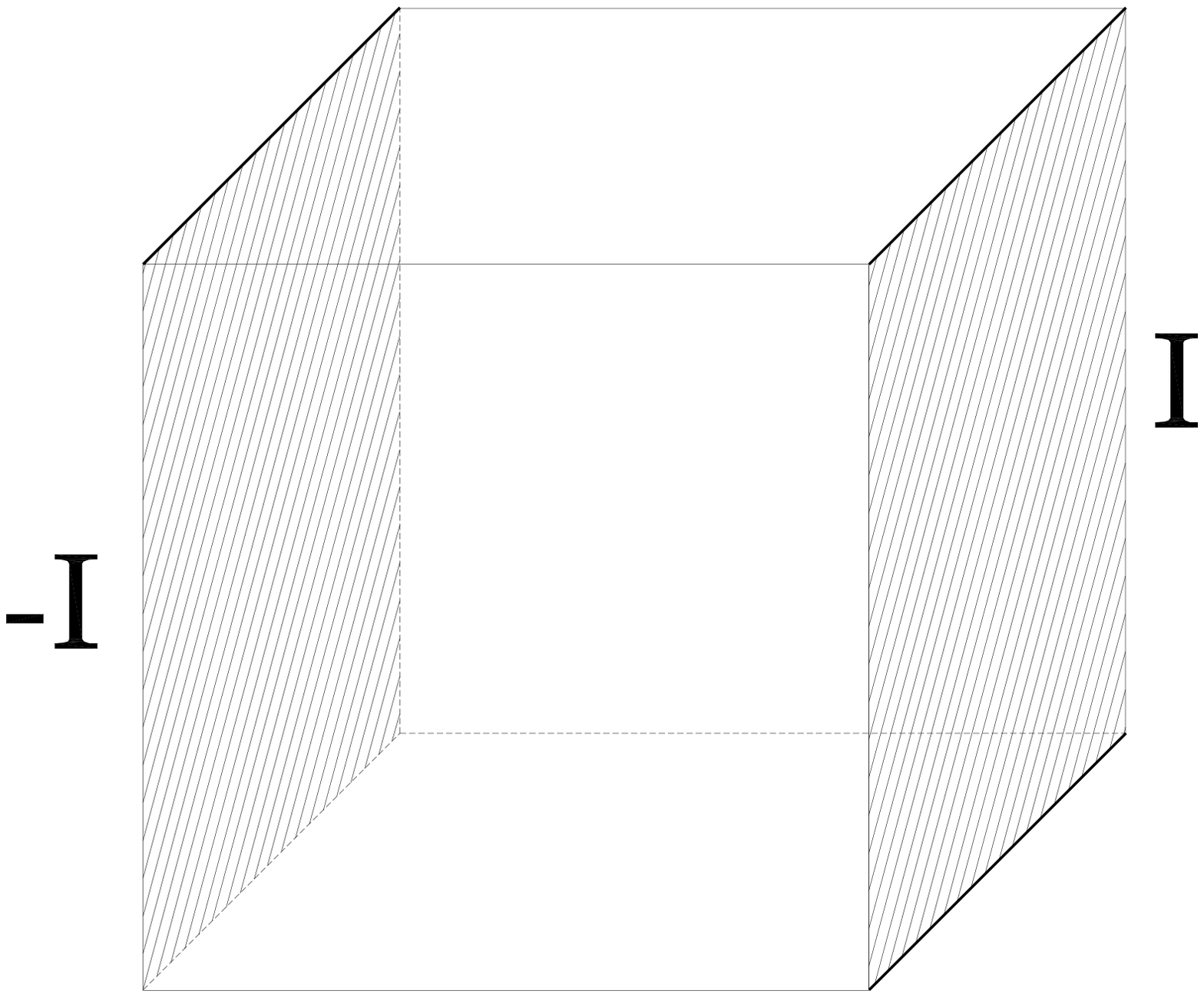}\\
     \centering{(a)}
   \end{minipage}
   \begin{minipage}{.30\textwidth}
     \centering
     \includegraphics[width=4cm]{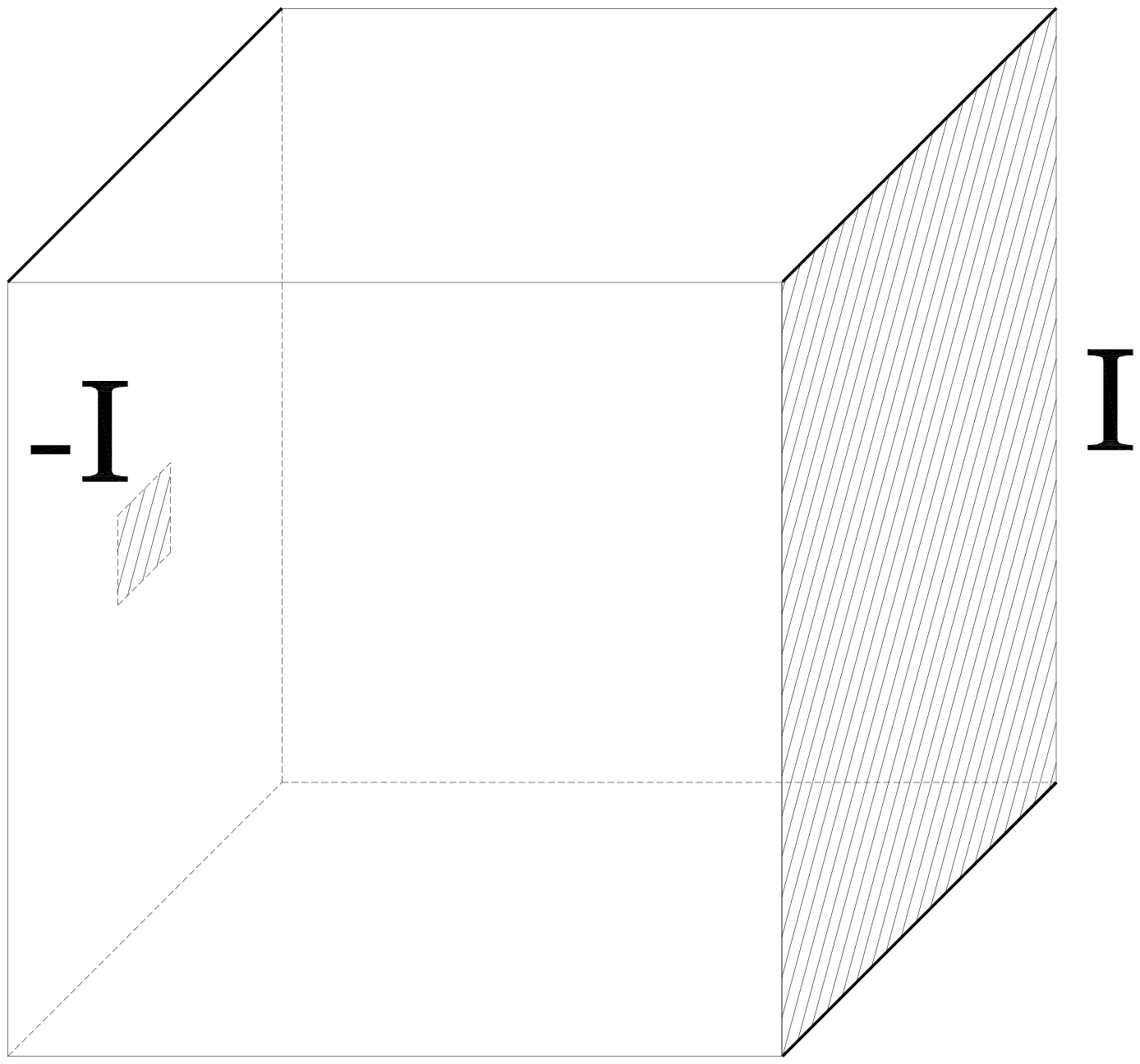}\\
     \centering{(b)}
   \end{minipage}
   \begin{minipage}{.30\textwidth}
     \centering
     \includegraphics[width=4cm]{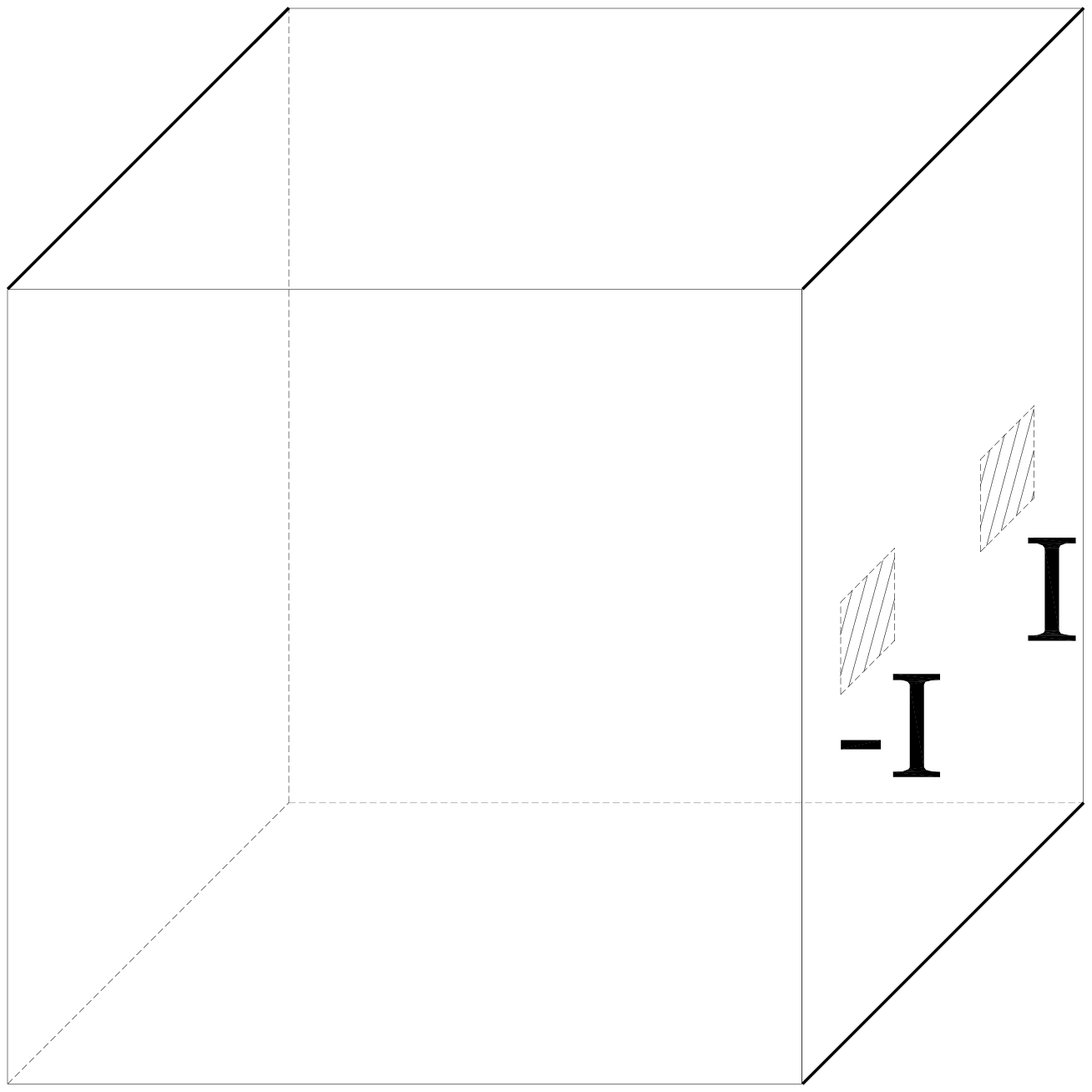}\\
     \centering{(c)}
   \end{minipage}
   \caption{Cubic conductor considered in 3--D numerical simulations for the physical  EIT model and location of the electrodes: test $T_1$ (a), test  $T_2$ (b) and test $T_3$ (c).} 
    \label{fig:EIT_prototype}
\end{figure}

The analysis has been restricted to the case of two electrodes located on the boundary of a cubic electrical conductor of side
$l$, see Figure \ref{fig:EIT_prototype}. The specimen has been discretized by a mesh of $17 \times 17 \times 17$ cubic HC finite
elements and the numerical experiments have been carried out on cubic inclusions only, with volume up to $6\%$ of the total volume and conductivity value $k=0.1$ or $k=10$. The surface impedance takes a constant value such that $\zeta=\frac{z \sigma}{l}=0.2$ on both electrodes, according to properties of human skin reported in literature, see, for instance, \cite{l:ssbs}.

In test $T_1$ of Figure \ref{fig:EIT_prototype}, the electrodes cover completely two opposite faces of the specimen, whereas in Test $T_2$ one electrode coincides with a face of $\partial \Omega$ and the other is a square, formed by one or nine surface finite elements, and it is located in central position of the opposite face.  Finally, in Test $T_3$, two electrodes are placed on the same face of the conductor $\Omega$ in a symmetric way respect to middle lines of the face. The electrodes are separated by three finite elements and their dimensions are equal to the element size.

\begin{figure}[h]
\begin{minipage}{0.49\textwidth}
    \includegraphics[width=6cm]{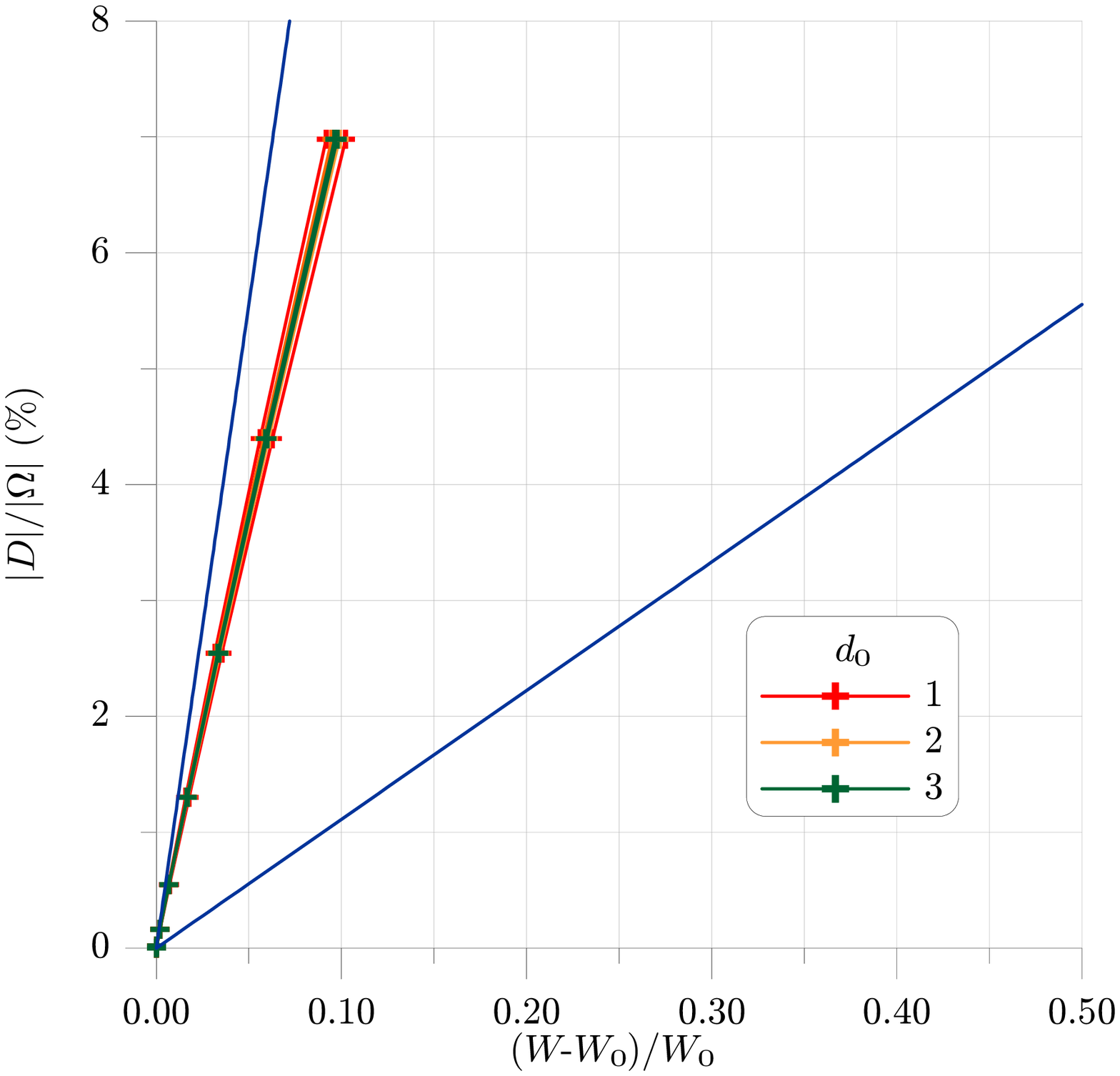}
    \centering{(a)}
\end{minipage}
\begin{minipage}{0.49\textwidth}
    \includegraphics[width=6cm]{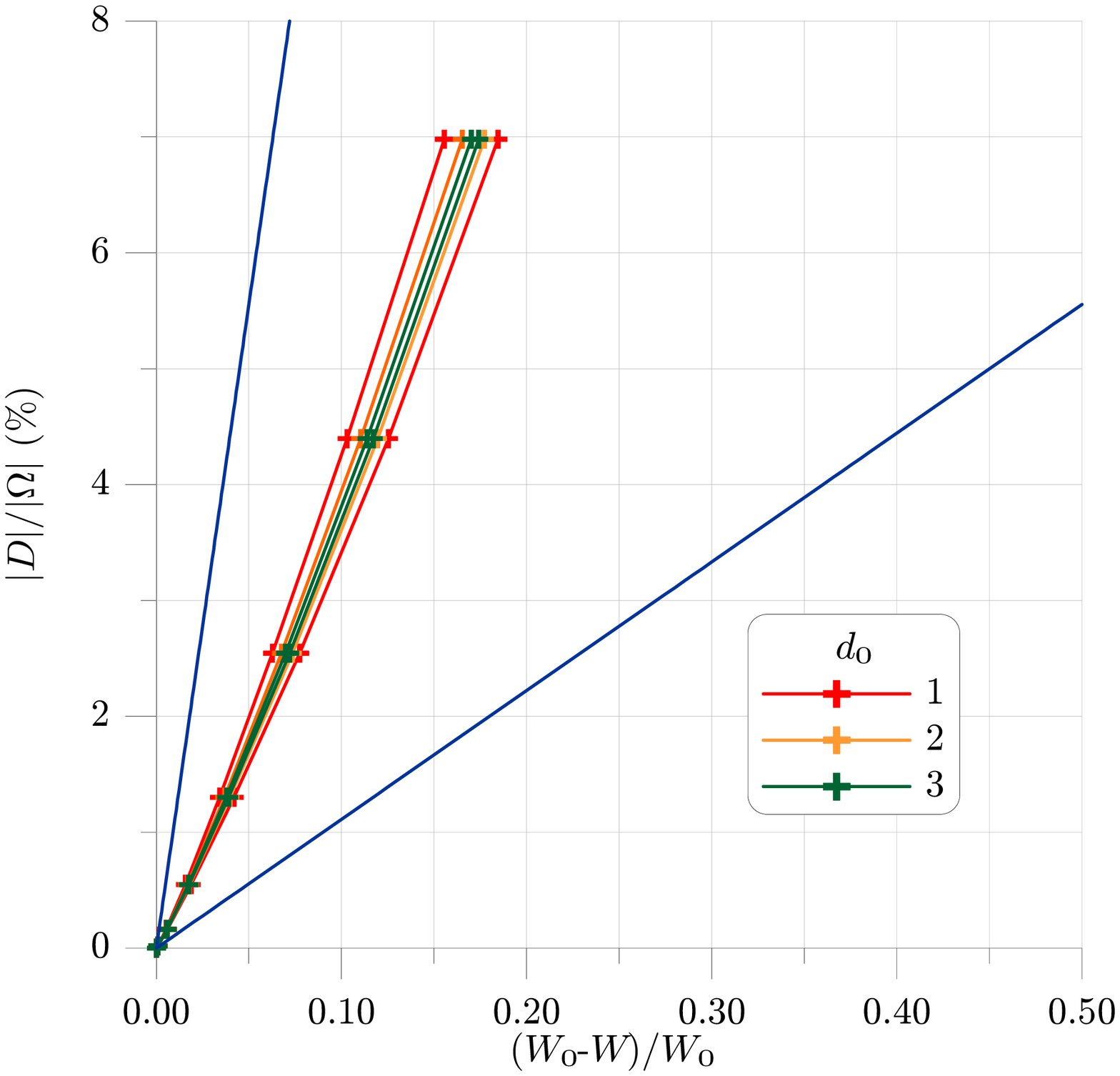}
    \centering{(b)}
\end{minipage}
\caption{Influence of $d_0$ for cubic inclusions in Test $T_1$ of
Figure \ref{fig:EIT_prototype}(a) ($17 \times 17 \times 17$ FE
mesh, $\zeta=0.2$ ): $k=0.1$ (a), $k=10$ (b).} \label{fig:T1_reg}
\end{figure}

The numerical results for Test $T_1$ are presented in Figure \ref{fig:T1_reg} for $k=0.1$ and $k=10$, respectively, and for
varying values of $d_0$. For both cases $k=0.1$ and $k=10$, the theoretical size estimates are given by
\begin{equation*}
   \frac{1}{9} \left ( \frac{l+2z}{l} \right )
   \frac{|W-W_0|}{W_0} \leq \frac{|D|}{|\Omega|}\leq \frac{10}{9}
   \left ( \frac{l+2z}{l} \right )
   \frac{|W-W_0|}{W_0}
\end{equation*}
and, again, they lead to a rather pessimistic evaluation of the
upper and lower bounds.

\begin{figure}[h]
\begin{minipage}{0.49\textwidth}
    \includegraphics[width=6cm]{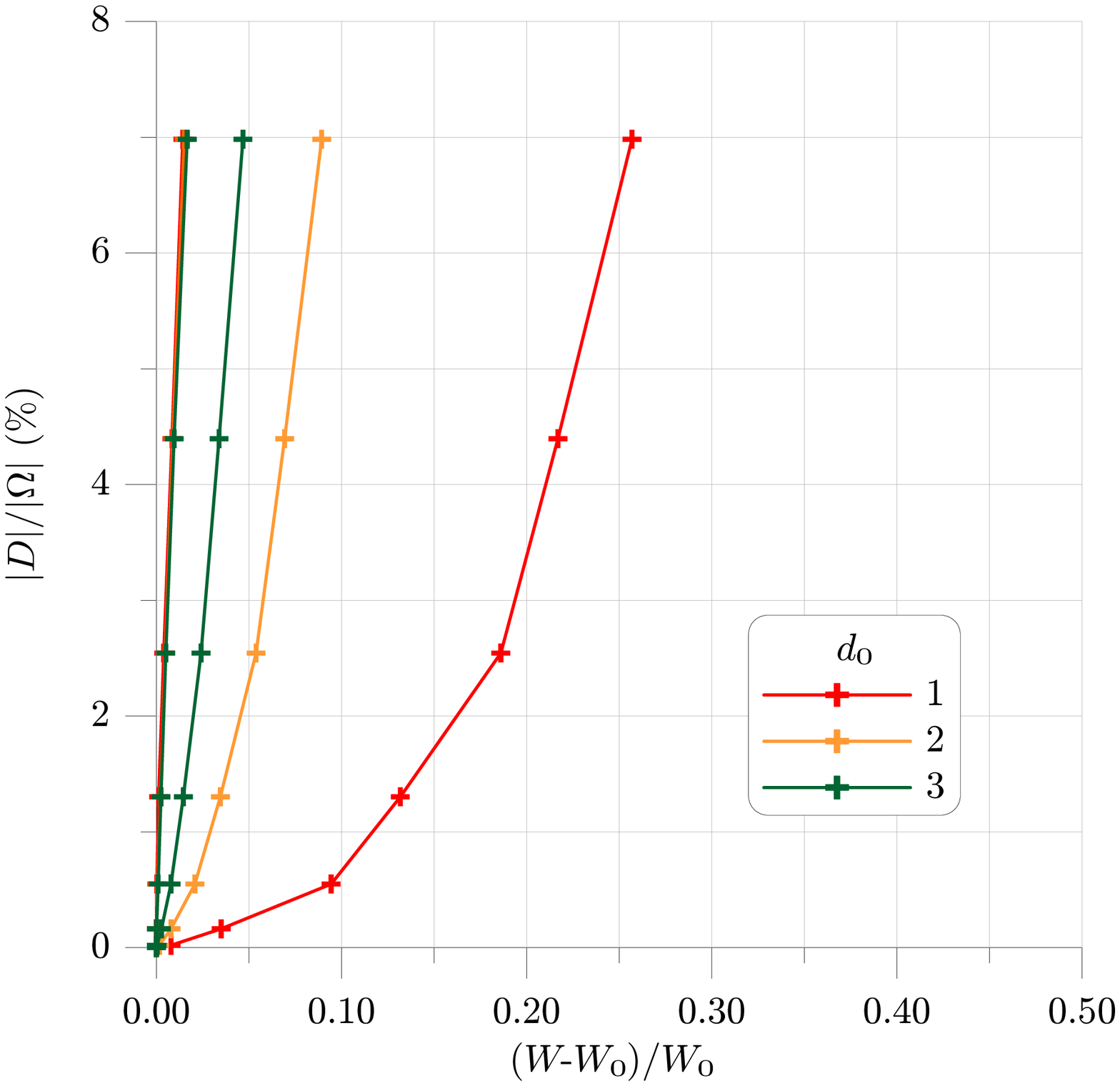}
     \centering{(a)}
\end{minipage}
\begin{minipage}{0.49\textwidth}
    \includegraphics[width=6cm]{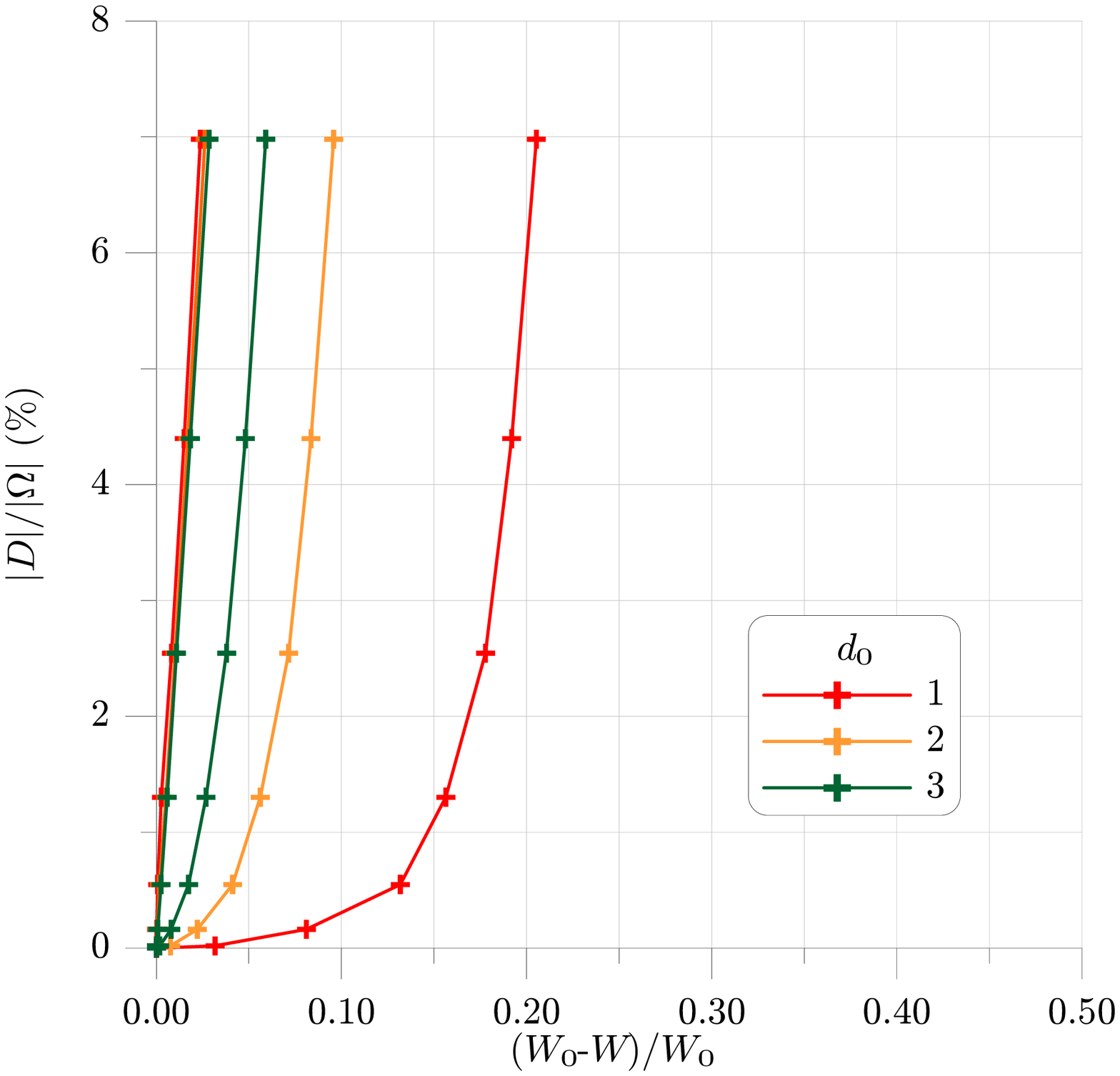}
    \centering{(b)}
\end{minipage}
\caption{Influence of $d_0$ for cubic inclusions in Test $T_2$ of
Figure \ref{fig:EIT_prototype}(b) ($17 \times 17 \times 17$ FE
mesh, $\zeta=0.2$, $1 \times 1$ FE electrode): $k=0.1$ (a), $k=10$
(b).} \label{fig:T2_reg_1}
\end{figure}

\begin{figure}[h]
\begin{minipage}{0.49\textwidth}
    \includegraphics[width=6cm]{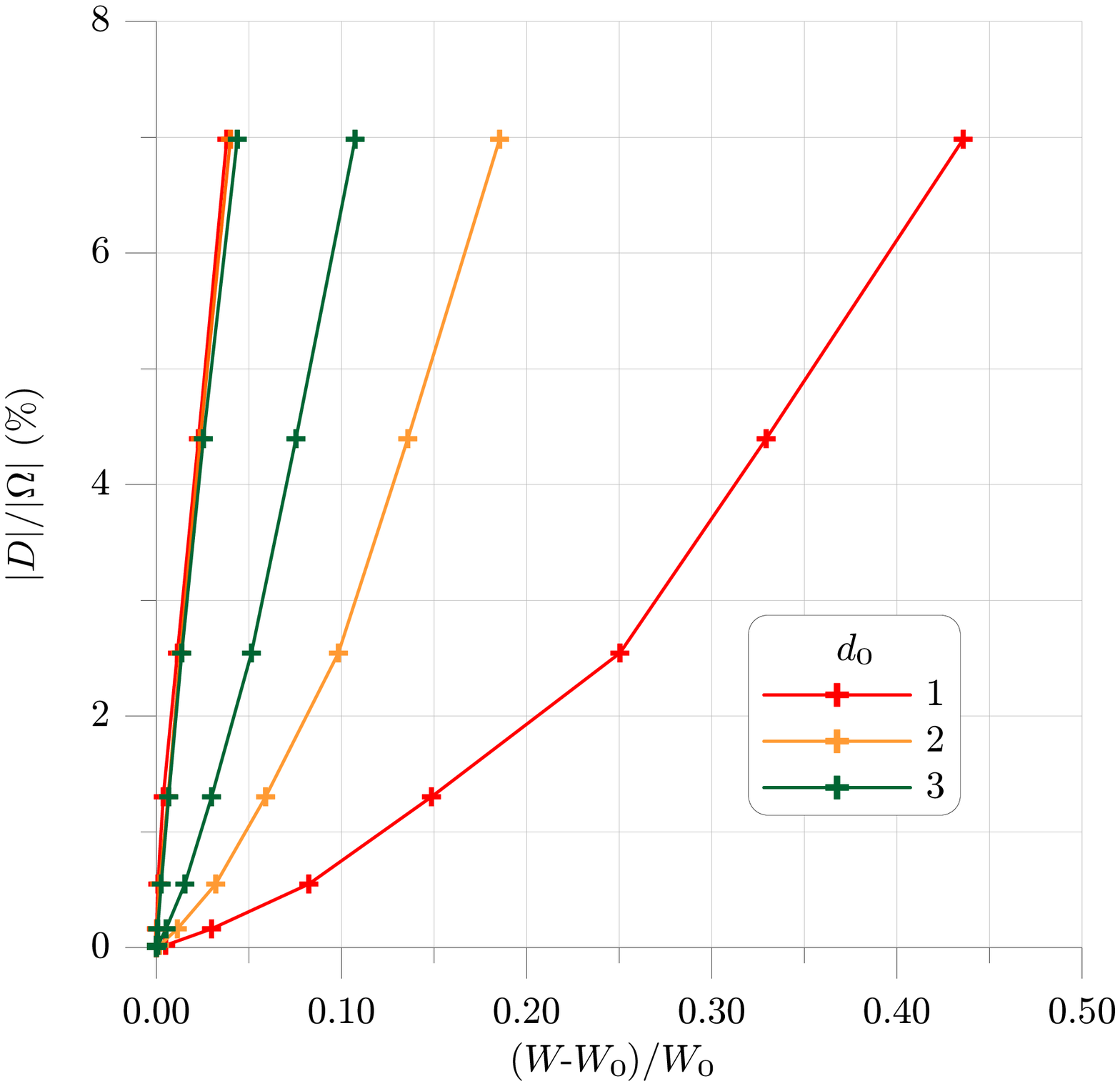}
    \centering{(a)}
\end{minipage}
\begin{minipage}{0.49\textwidth}
    \includegraphics[width=6cm]{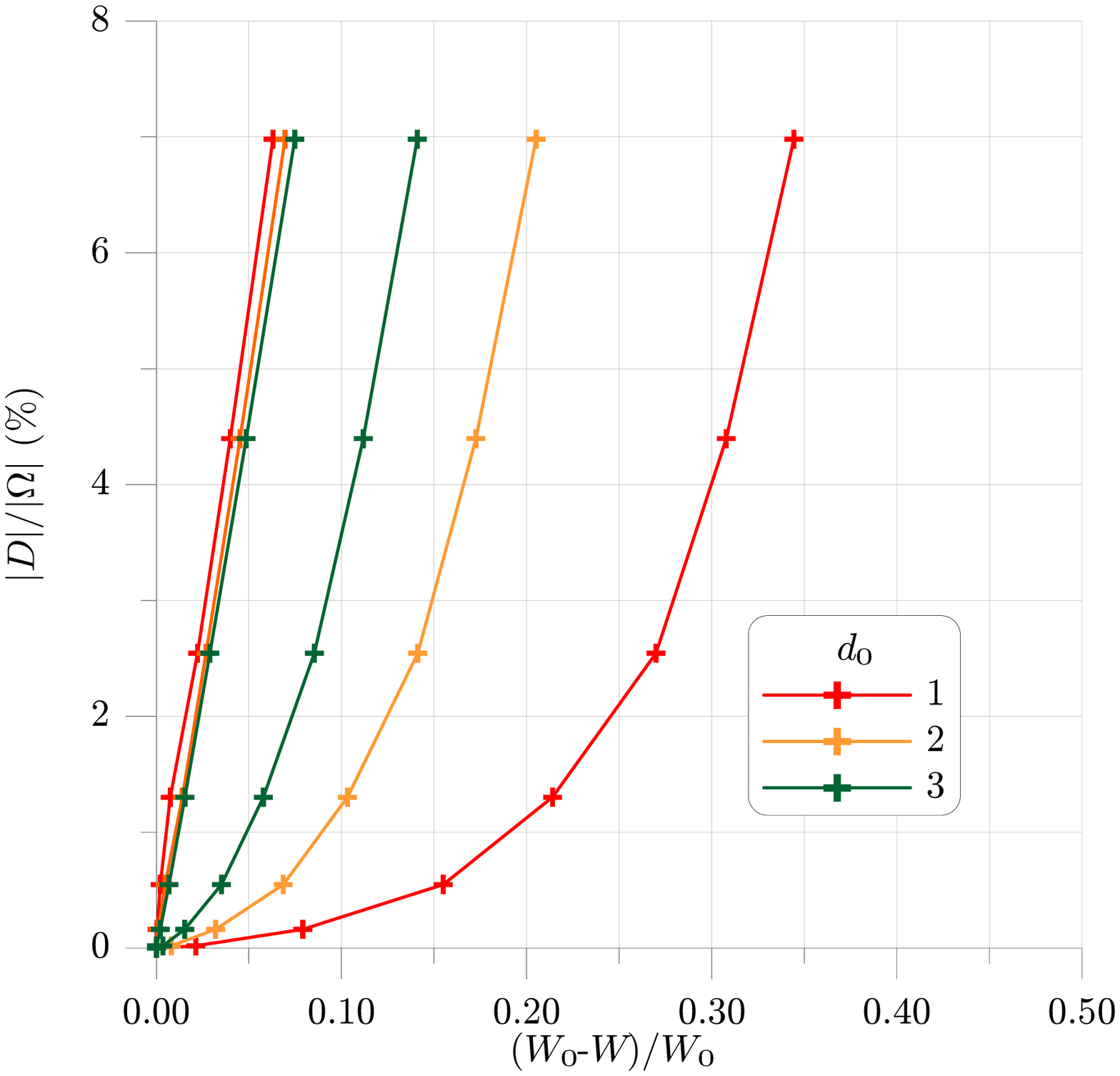}
    \centering{(b)}
\end{minipage}
\caption{Influence of $d_0$ for cubic inclusions in Test $T_2$ of
Figure \ref{fig:EIT_prototype}(b) ($17 \times 17 \times 17$ FE
mesh, $\zeta=0.2$, $3 \times 3$ FE electrode): $k=0.1$ (a), $k=10$
(b).} \label{fig:T22_reg_3}
\end{figure}

Concerning Test $T_2$, Figure \ref{fig:T2_reg_1} shows the results
when the small electrode coincides with one surface finite
element, whereas Figure \ref{fig:T22_reg_3} refers to the case of
a $3 \times 3$ finite elements electrode. One can notice that in
all the four cases, the upper bound is not really influenced by
the value of $d_0$. Moreover, the inaccuracy in determining the
lower bound of the angular sector, is probably due to the fact
that the present analysis is restricted to the special class of
cubic inclusions.

\begin{figure}[h]
\begin{minipage}{0.49\textwidth}
    \includegraphics[width=6cm]{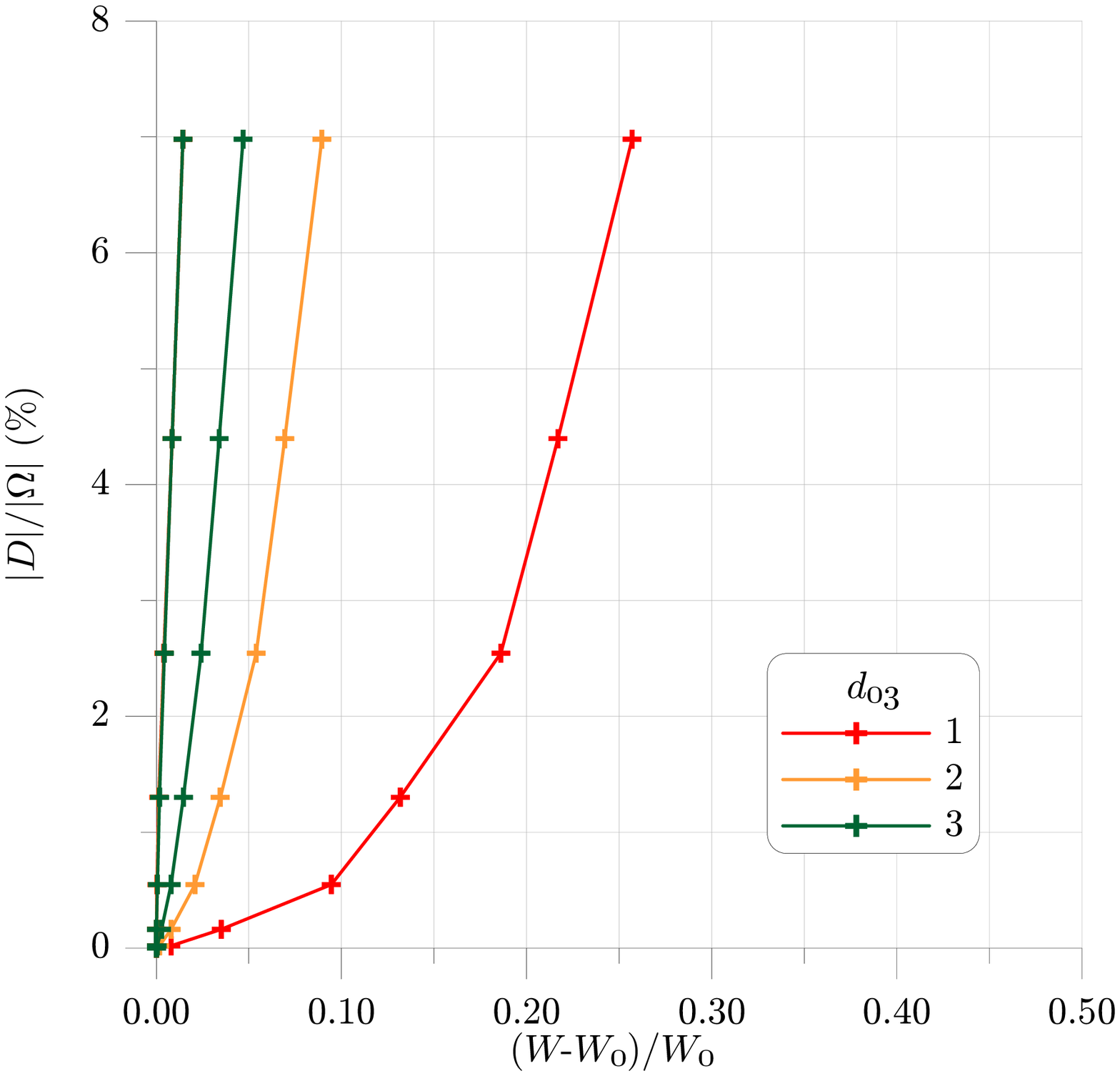}
     \centering{(a)}
\end{minipage}
\begin{minipage}{0.49\textwidth}
    \includegraphics[width=6cm]{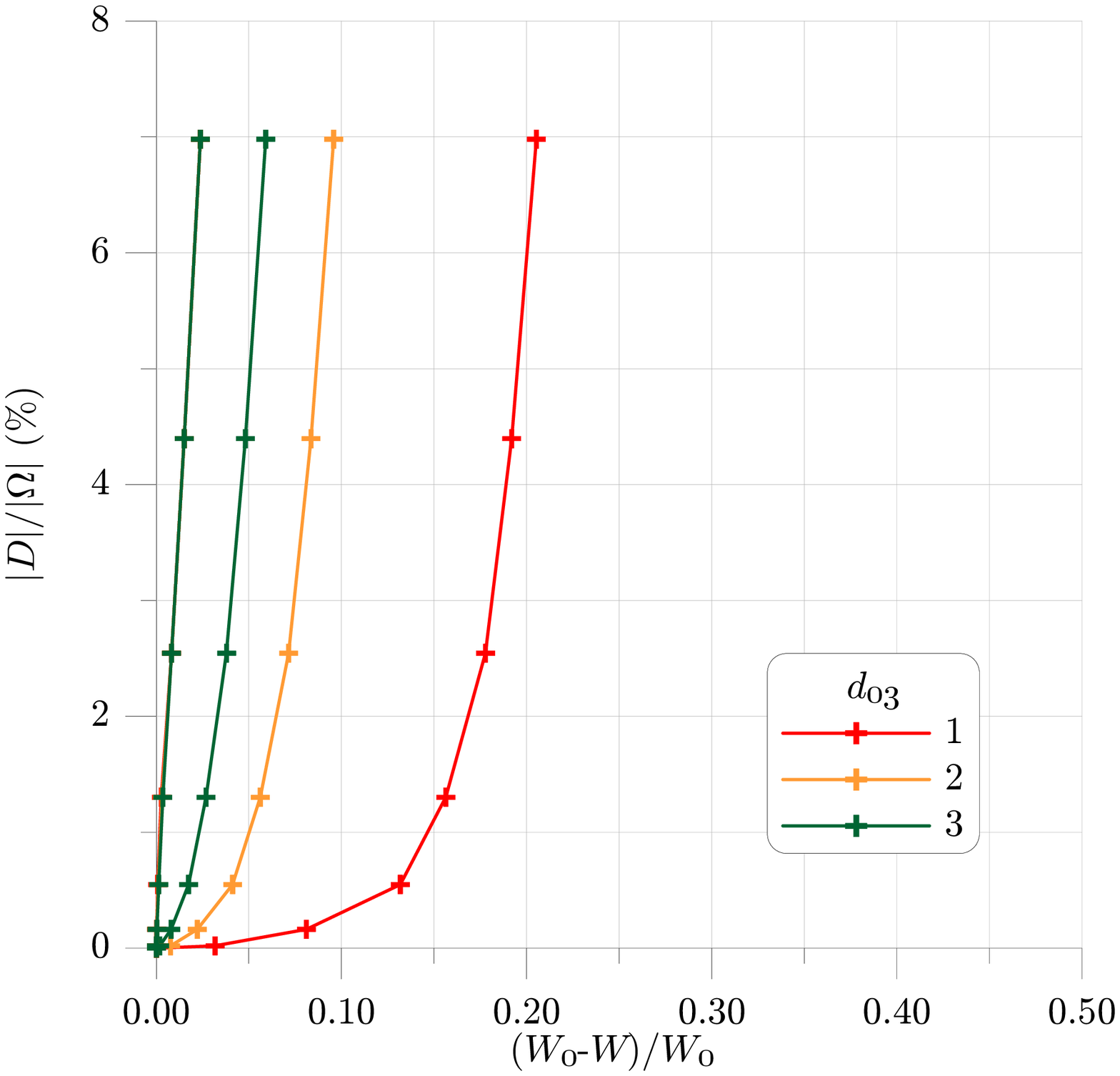}
    \centering{(b)}
\end{minipage}
\caption{Influence of $d_{03}$ for cubic inclusions in Test $T_2$
of Figure \ref{fig:EIT_prototype}(b) ($17 \times 17 \times 17$ FE
mesh, $\zeta=0.2$, $1 \times 1$ FE electrode): $k=0.1$ (a), $k=10$
(b).} \label{fig:T2_rvar}
\end{figure}

A comparison between Figure \ref{fig:T2_reg_1} and Figure
\ref{fig:T22_reg_3} suggests that better upper bounds can be
obtained by enlarging the size of the small electrode. Moreover,
from Figures \ref{fig:T2_reg_1} and \ref{fig:T22_reg_3} it appears
clearly that the lower bound significantly improves as the
distance $d_0$ between the inclusion $D$ and the boundary of
$\Omega$ increases. This property has been further investigated by
increasing only the distance $d_{03}$ of the inclusion $D$ from
the face of the conductor containing the small electrode. Figure
\ref{fig:T2_rvar} shows the results of simulations in the case of
a single finite element electrode and a comparison with Figure
\ref{fig:T2_reg_1} suggests that the improvement of the lower
bound is mainly due to the greater distance from the electrode.


\begin{figure}[h]
\begin{minipage}{0.49\textwidth}
    \includegraphics[width=6cm]{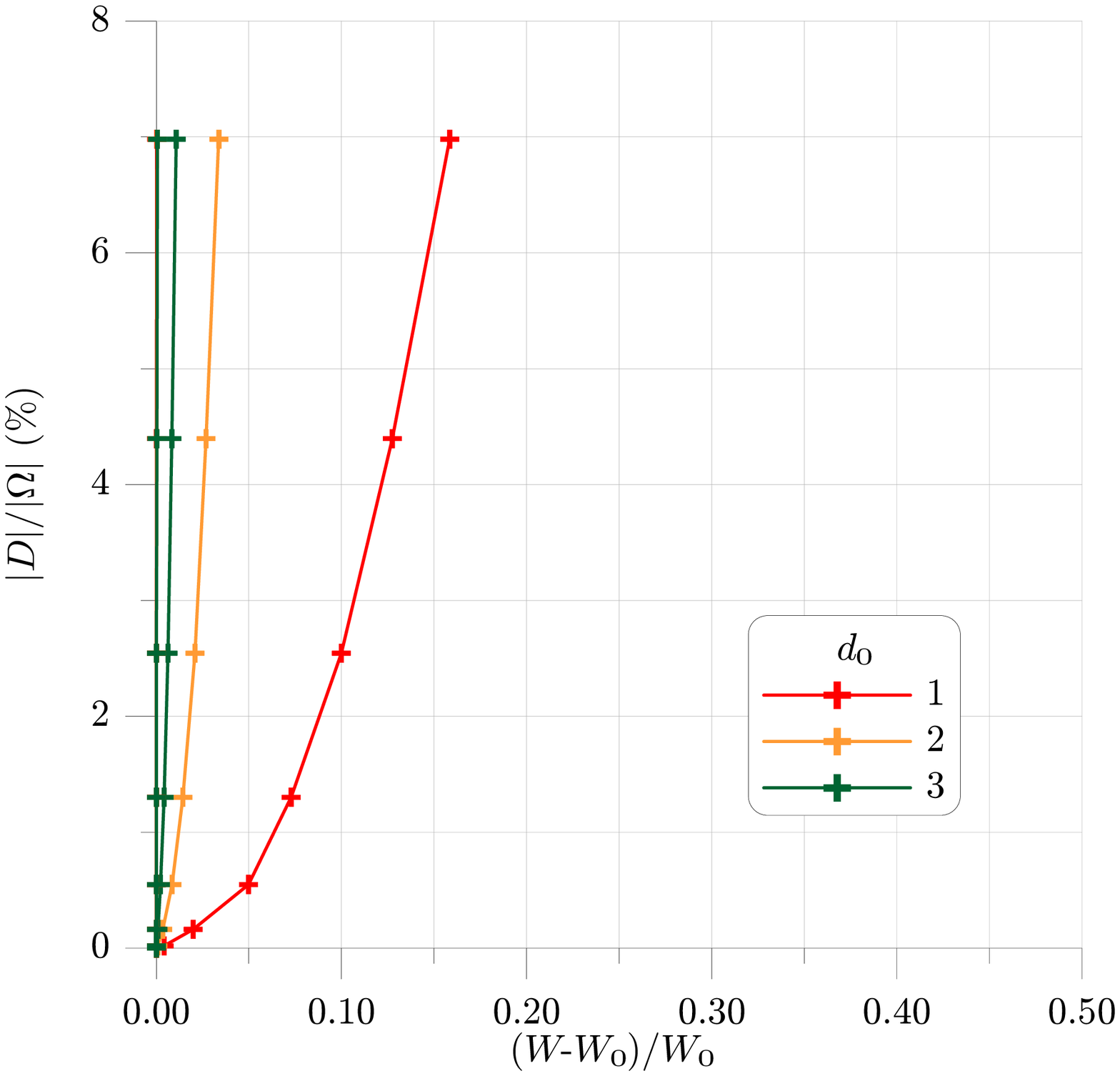}
     \centering{(a)}
\end{minipage}
\begin{minipage}{0.49\textwidth}
    \includegraphics[width=6cm]{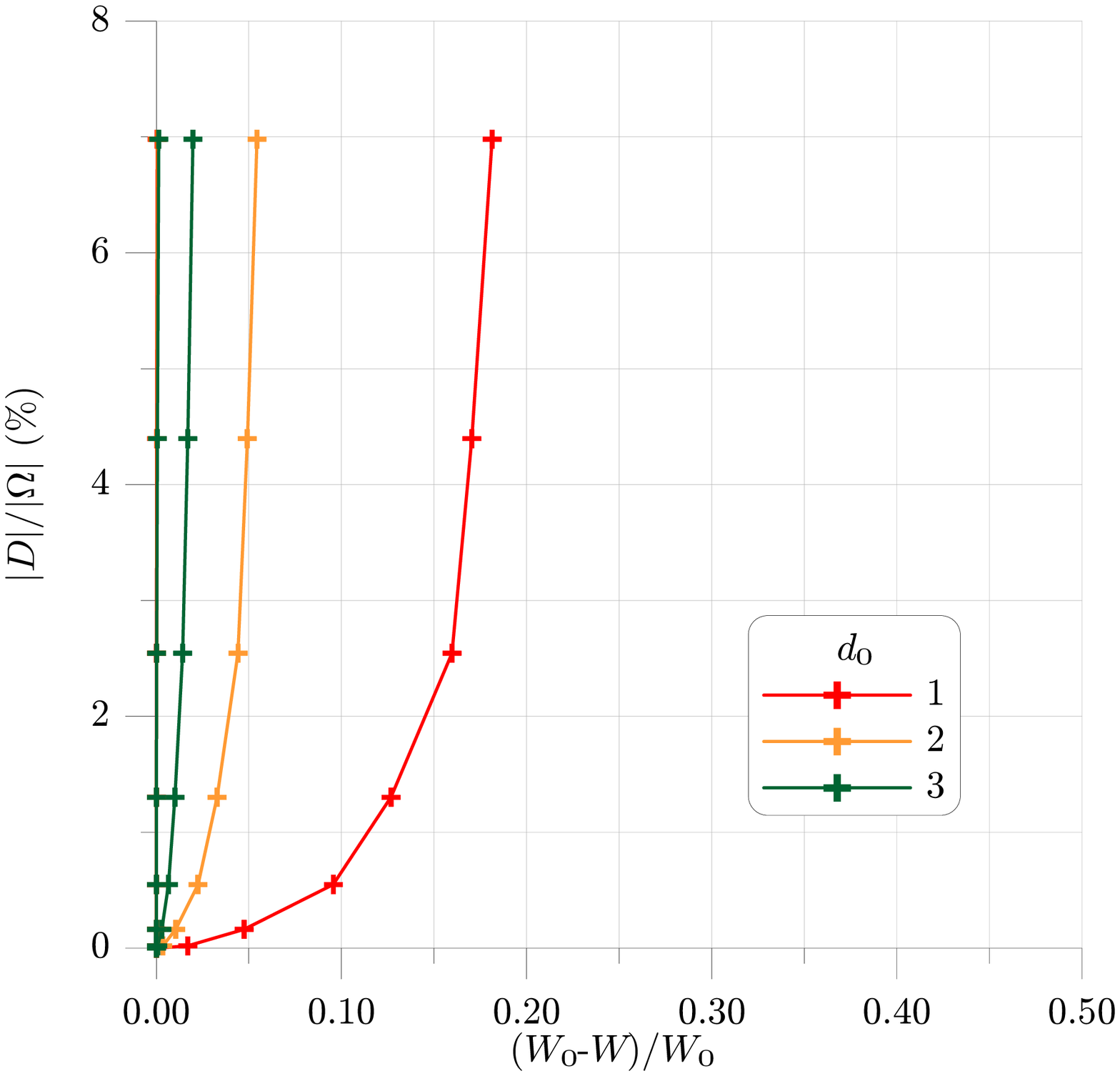}
    \centering{(b)}
\end{minipage}
\caption{Influence of $d_0$ for cubic inclusions in Test $T_3$ of Figure \ref{fig:EIT_prototype}(c) ($17 \times 17 \times 17$ FE
mesh, $\zeta=0.2$, $1 \times 1$ FE electrode): $k=0.1$ (a), $k=10$
(b).} \label{fig:T3_reg}
\end{figure}

Finally, the results of the numerical simulations for Test $T_3$ are presented in Figure \ref{fig:T3_reg}. In this case, the lower bound improves as the distance $d_0$ between the inclusion $D$ and the boundary of $\Omega$ increases, whereas the upper bound is indistinguishable from the vertical axis.

\section{Conclusions} \label{sec:conclusions}

We have tested by numerical simulations the approach of \emph{size
estimates} for EIT. We could perform experiments in the 2--D
setting with a large varieties of shapes of inclusions and we
found quite satisfactory bounds, which in some cases are markedly
better than those derived theoretically.

In the 3--D case, we had to limit the variety of shapes of the
test inclusions since the growth of their degree of freedom
conflicts with the limitations on computer time. We showed that
good volume bounds hold when the boundary data $\varphi$ is
\emph{well-behaved} in terms of its frequency, whereas they
rapidly deteriorate as the frequency increases.

For the complete EIT model we have also made tests in a 3--D
setting and compared the bounds in terms of the size of the
electrodes, their relative distance and their a-priori assumed
distance from the inclusion $D$. We have shown that we obtain good
bounds when the electrodes are not too small and when $D$ is
sufficiently away from them.

%
\bibliographystyle{alpha}
%
%

\end{document}